\documentclass{article}
\usepackage[T2A]{fontenc}
\usepackage[cp1251]{inputenc}
\usepackage[english,russian]{babel}
\usepackage[tbtags]{amsmath}
\usepackage{amsfonts,amssymb,mathrsfs,amscd}
\numberwithin{equation}{section}
\newtheorem{remark}{Remark}[section]
\newtheorem{lemma}{Lemma}[section]
\newtheorem{theorem}{Theorem}[section]

\newtheorem{conclusion}{Conclusion}

\newcommand{\col}{\mathop{\rm col}}

\begin{document}
\begin{Large}
\thispagestyle{empty}
\begin{center}
{\bf Inverse scattering problem for a cubic string having the shape of a step\\
\vspace{5mm}
V. A. Zolotarev}\\

B. Verkin Institute for Low Temperature Physics and Engineering
of the National Academy of Sciences of Ukraine\\
47 Nauky Ave., Kharkiv, 61103, Ukraine

\end{center}
\vspace{5mm}

{\small {\bf Abstract.} Problem for a cubic string having the shape of a step is studied. It is proved that here we have two scattering problems: the direct one and its dual. The direct problem describes scattering of the waves coming from ``$+\infty$'', and its dual, correspondingly, of the waves from ``$-\infty$''. The main systems of linear singular equations for both problems are obtained. It is shown how to unambiguously restore potential from the solutions to these equation systems.}
\vspace{5mm}


\vspace{5mm}

\begin{center}
{\bf Introduction}
\end{center}
\vspace{5mm}

Inverse scattering problem is one of the most interesting and intriguing sections of mathematical physics and functional analysis. It emerged in the classic works of V. A. Marchenko \cite{1} and M. G. Krein (\cite{1} -- \cite{4}). For the Sturm -- Liouville operator, procedure of potential restoration from the scattering data was found, this procedure is based on the Marchenko equation \cite{1}.

Method of inverse scattering problem is applied wondrously and effectively in the theory of non-linear partial differential equations \cite{5}. A non-linear equation written in the Lax form, $L_t=[L,A]$ ($L$, $A$ are linear differential operators), as a rule, led to the fact that $L$ is a Surm -- Liouville operator with the solution to a non-linear equation as a potential. Next, knowing the evolution by $t$ of the scattering data, the non-linear equation was solved using the inverse scattering problem \cite{5}.

For several equations (Degasperis -- Processi), search for the Lax pair $L$, $A$  leads to the third-order operator $L$ \cite{6} -- \cite{10}. The main equation here is the equation of a cubic string
\begin{equation}
y'''(x)=m(x)\lambda^3y(x),\label{eq0.1}
\end{equation}
its name appears analogously to the classic studies for a usual string \cite{9,10}. Spectral problem for the cubic string \eqref{eq0.1} is studied in the works \cite{12,13}.

The present work studies scattering problem for the equation of the form
\begin{equation}
iy'''(x)=m(x)\lambda^3y(x)\quad(x\in\mathbb{R},\lambda\in\mathbb{C})\label{eq0.2}
\end{equation}
where potential $m(x)$ is a real function such that

(i) the limits $\lim\limits_{x\rightarrow\pm\infty}m(x)=m_{\pm}$ exist ($0<m_{\pm}<\infty$);

(ii) there exists such $a>0$ that
\begin{equation}
\int\limits_{\mathbb{R}_\pm}|m(x)-m_\pm|^2e^{2ax}dx<\infty.\label{eq0.3}
\end{equation}
Difficulty of the studied scattering problem consists in the fact that the potential $m(x)$ enters equation \eqref{eq0.1} in a multiplicative way and not additively, as usual \cite{14, 15}.

Operator corresponding to equation \eqref{eq0.2} is not-self-adjoint, therefore here we study two scattering problems corresponding to the waves incident from ``$+\infty$'' and from ``$-\infty$''. These problems are dual to each other.

Section 1 proves existence of the Jost solutions at ``$+\infty$'' behaving asymptotically as $\left\{e^{i\lambda n_+\zeta_kx}\right\}_0^2$ where $\{\zeta_k\}_0^2$ are the cubic roots of the identity and $n_+^3=m_+$. In the same way, it is established existence of the Jost solution at ``$-\infty$'' with asymptotic $\{e^{i\lambda n_-\zeta_kv}\}_0^2$ ($n_-^3=m_-$). In both cases, analytic properties relative to $\lambda$ of such Jost solutions are described.

Section 2 studies the direct scattering problem for the waves incident from ``$+\infty$''. Properties of the transition matrix are established. It is shown that zeros of the holomorphic function that correspond to the bound states have multiplicity 2. Jump problems on the rays originating at origin are constructed. The main closed system of linear singular equations for the inverse problem (an analogue to the Marchenko equation) is obtained. Scattering coefficients and spectrum points connected with bound states are the independent parameters of this system. Knowing solution to this main equation system, using the simple formula, we restore the potential $m(x)$ on the right half-axis $x\in\mathbb{R}_+$.

Section 3 deals with the dual scattering problem corresponding to the waves incident from ``$-\infty$''. Similarly to Section 2, for scattering coefficients, we obtain relation corresponding to the unitarity of the problem. Moreover, the relation establishing reciprocity between direct and dual scattering systems is found. This section also derives the main closed system of linear singular integral equations for the dual scattering problem. From the solution to this system, using a simple formula, we find potential $m(x)$ on the left half-axis $x\in\mathbb{R}_-$.

So, from the solutions to the two main equation systems, potential $m(x)$ is restored on the entire axis ($x\in\mathbb{R}$). To obtain the main equation system, Sections 3 and 4 use the Riemann boundary value problem on the system of rays outgoing from origin.

In conclusion, we notice that results of the Sections 1 -- 3 are based upon and use methods and techniques of the solution of inverse scattering problems for third-order operators obtained in \cite{14,15}.

\section{Preliminary information and Jost solutions}\label{s1}

{\bf 1.1} Equation
\begin{equation}
iy'''(x)=\lambda^3y(x)\quad(x\in\mathbb{R},\lambda\in\mathbb{C})\label{eq1.1}
\end{equation}
has three independent solutions $\{e^{i\lambda\zeta_kx}\}_0^2$ where $\{\zeta_k\}_0^2$ are the roots of the cubic equation $\zeta^3=1$, and
\begin{equation}
\zeta_0=1,\quad\zeta_1=\frac12(-1+i\sqrt3),\quad\zeta_2=\frac12(-1-i\sqrt3).\label{eq1.2}
\end{equation}
Another system of fundamental solutions $\{s_k(i\lambda x)\}_0^2$ \cite{14,15} of equation \eqref{eq1.1} where
\begin{equation}
s_k(z)=\frac13\left(\frac1{\zeta_0^k}e^{\zeta_0z}+\frac1{\zeta_1^k}e^{\zeta_1z}+\frac1{\zeta_2^k}e^{\zeta_2z}\right)\label{eq1.3}
\end{equation}
plays an important role. Properties of the functions $\{s_k(z)\}$ are given in \cite{14,15}, these functions are analogues of sines and cosines for equation \eqref{eq1.1}. Among all these properties, we list the following \cite{14,15}:

(i) $s_k(\zeta_1z)=\zeta_1^ks_k(z)$ ($k$-evenness);

(ii) the Euler formula
$$e^{z\zeta_k}=s_0(z)+\zeta_ks_1(z)+\zeta_k^2s_2(z);$$

(iii) the main identity
\begin{equation}
s_0^3(z)+s_1^3(z)+s_2^3(z)-3s_0(z)s_1(z)s_2(z)=1;\label{eq1.4}
\end{equation}

(iv) the functions $\{s_k(z)\}_0^2$ are solutions to the differential equation $y'''(z)=y(z)$ and satisfy the initial conditions
$$s_0(0)=0;\quad s'_0(0)=0;\quad s''_0(0)=0;$$
$$s_1(0)=0;\quad s'_1(0)=1;\quad s''_1(0)=0;$$
$$s_2(0)=0;\quad s'_2(0)=0;\quad s''_2(0)=0.$$

Consider the Cauchy problem
\begin{equation}
iy'''(x)=\lambda^3y(x)=f(x);\quad y(0)=y_0,\,y'(0)=y_1,\,y''(0)=y_2\,(\lambda\in\mathbb{C},x\in\mathbb{R}_+).\label{eq1.5}
\end{equation}
Using (iv) \eqref{eq1.4}, we find the solution to the homogenous ($f(x)\equiv0$) problem \eqref{eq1.5},
\begin{equation}
y_0(\lambda,x)=y_0s_0(i\lambda x)+y_1\frac{s_1(i\lambda x)}{i\lambda}+y_2\frac{s_2(i\lambda x)}{i\lambda x}\label{eq1.6}
\end{equation}
General solution to the Cauchy problem \eqref{eq1.5} equals
\begin{equation}
y(\lambda,x)=y_0s_0(i\lambda x)+y_1\frac{s_1(i\lambda x)}{i\lambda}+y_2\frac{s_2(i\lambda x)}{(i\lambda)^2}+i\int\limits_0^x\frac{s_2(i\lambda x)}{(i\lambda)^2}f(t)dt.\label{eq1.7}
\end{equation}
Unit vectors $\{\zeta_k\}_0^2$ \eqref{eq1.2} are the direction vectors of three straight lines in $\mathbb{C}$,
\begin{equation}
L_{\zeta_k}=\{x\zeta_k: x\in\mathbb{R}\}\quad(0\leq k\leq2).\label{eq1.8}
\end{equation}
By $\{l_{\zeta_k}\}_0^2$, we denote rays in the direction $\zeta_k$ outgoing from origin, and by $\{\widehat{l}_{\zeta_k}\}_0^2$, rays in the direction $\zeta_k$, coming to origin,
\begin{equation}
l_{\zeta_k}=\{x\zeta_k:x\in\mathbb{R}_+\},\quad\widehat{l}_{\zeta_k}=L_{\zeta_k}\setminus l_{\zeta_k}=\{x\zeta_k:x\in\mathbb{R}_-\}.\label{eq1.9}
\end{equation}
The straight lines $\{L_{\zeta_k}\}_0^3$ \eqref{eq1.8} divide the $\mathbb{C}$ plane into six sectors:
\begin{equation}
S_p=\left\{z\in\mathbb{C}:p\frac{2\pi}6<\arg z<(p+1)\frac{2\pi}6\right\}\quad(p=0,1,...,5).\label{eq1.10}
\end {equation}
\vspace{5mm}

{\bf 1.2} Consider equation of the cubic string \cite{9,10,13},
\begin{equation}
iy'''(x)=m(x)\lambda^3y(x)\quad(x\in\mathbb{R},\lambda\in\mathbb{C})\label{eq1.11}
\end{equation}
where function $m(x)$ is real and satisfies the conditions

(i) the limits $\lim\limits_{x\rightarrow\pm\infty}m(x)=m_\pm$ exist ($|m_\pm|<\infty$, $m_\pm=n_\pm^3$, $n_\pm\in\mathbb{R}$);

(ii) for some $a>0$,
\begin{equation}
\int\limits_{\mathbb{R}_\pm}|m(x)-m_\pm|^2e^{2a|x|}dx<\infty.\label{eq1.12}
\end{equation}

\begin{remark}\label{r1.1}
Hereinafter, we assume that $m_+$ and $m_-$ are of the same sign, $m_+>0$ and $m_->0$. Condition (ii) \eqref{eq1.12} implies that $|m(t)-m_\pm|e^{b|x|}\in L^1(\mathbb{R}_\pm)\cap L^2(\mathbb{R}_\pm)$, for all $b$ such that $b<a$.
\end{remark}

Equation \eqref{eq1.11} has a natural operator interpretation. Let $m(x)\geq m_0>0$, define the weight space $L^2(m,\mathbb{R})$ formed by the complex-valued functions $f(x)$ such that
\begin{equation}
\int\limits_\mathbb{R}|f(x)|^2m(x)dx<\infty\label{eq1.13}
\end{equation}
and consider the self-adjoint operator
\begin{equation}
(\mathcal{L}_my)(x)=\frac1{m(x)}iy'''(x)\label{eq1.14}
\end{equation}
in this space, with natural domain $\mathfrak{D}(\mathcal{L}_n)=W_2^3(m,\mathbb{R})$.
In terms of $\mathcal{L}_m$ \eqref{eq1.14}, equation \eqref{eq1.11} becomes
\begin{equation}
\mathcal{L}_my=\lambda^3y.\label{eq1.15}
\end{equation}

\begin{remark}\label{r1.2}
Self-adjointness of $\mathcal{L}_m$ \eqref{eq1.14} and \eqref{eq1.15} imply that $\lambda^3$ is real. Therefore, hereinafter we assume that $\lambda^3$ in equation \eqref{eq1.11} is real.
\end{remark}

Condition (i) \eqref{eq1.12} implies that, for $x\rightarrow\pm\infty$, equation \eqref{eq1.11} transforms into $iy'''=(\lambda n_+)^3y$, equivalent to \eqref{eq1.1}. Therefore, it is natural to define the Jost solutions $\{v_k(\lambda,x)\}_0^2$ and $\{u_k(\lambda,x)\}+0^2$ as the solutions to equation \eqref{eq1.11} which have the following asymptotic behavior:
\begin{equation}
\begin{array}{ccc}
(a)\,v_k(\lambda,x)\rightarrow e^{iz_+\zeta_kx}&(x\rightarrow\infty,z_+=\lambda n_+,0\leq k\leq2);\\
(b)\,u_k(\lambda,x)\rightarrow e^{iz_-\zeta_kx}&(x\rightarrow-\infty,z_-=\lambda n_-,0\leq k\leq2).
\end{array}\label{eq1.16}
\end{equation}
The boundary value problem
\begin{equation}
\left\{
\begin{array}{lll}
iy'''(x)+q(x)y(x)=z^3y(x);\\
y(x)\rightarrow e^{iz\zeta x}\,(x\rightarrow\infty),
\end{array}\right.\label{eq1.17}
\end{equation}
as is known \cite{14,15} (see also \eqref{eq1.7}), is equivalent to the integral equation
\begin{equation}
y(\lambda,x)=e^{i\lambda\zeta x}-i\int\limits_x^\infty\frac{s_2(iz(x-t))}{(iz)^2}q(t)y(t)dt.\label{eq1.18}
\end{equation}
Rewrite equation \eqref{eq1.11} as
$$iy'''-\lambda^3(m(x)-m_+)y=z_+^3y.$$
Hence, using \eqref{eq1.18} for $\{v_k(\lambda,x)\}_0^2$, solutions to the boundary value problem \eqref{eq1.11} (a) \eqref{eq1.16} we obtain the integral equation
\begin{equation}
v_k(\lambda,x)=e^{iz_+\zeta_kx}-iz_+\int\limits_x^\infty s_2(iz_+(x-t))\cdot\left(\frac{m(t)}{m_+}-1\right)v_k(\lambda,t)dt\quad(0\leq k\leq2)\label{eq1.19}
\end{equation}
($z_+=\lambda n_+$). In $L^2(\mathbb{R})$, define the family of Volterra operators
$$(K_zf)(x)\stackrel{\rm def}{=}\int\limits_x^\infty K_1(z,x,t)\left(\frac{m(t)}{m_t}-1\right)f(t)dt,$$
where
\begin{equation}
K_1(z_+,x,t)=s_2(iz_+(x-t)).\label{eq1.20}
\end{equation}
Equation \eqref{eq1.19}, in terms of $K_{z_+}$ becomes
$$(I\cdot+izK_z)v_k(\lambda,x)=e^{iz_+\zeta_kx}\quad(0\leq k\leq2),$$
and thus
\begin{equation}
v_k(\lambda,x)=\sum\limits_{n=0}^\infty(-iz_+)^nK_{z_+}^ne^{iz_+\zeta_kx}.\label{eq1.21}
\end{equation}
Operators $K_{z_+}^n$ are also Volterra,
$$(K_{z_+}^nf)(x)=\int\limits_x^\infty K_n(z_+,x,t)\left(\frac{m(t)}{m_+}-1\right)f(t)dt,$$
and their kernels $K_n(z_+,x,t)$ satisfy the recurrent relations
\begin{equation}
K_{n+1}(z_+,x,t)=\int\limits_x^tK_n(z_+,x,s)\left(\frac{m(s)}{m_+}-1\right)K_1(z_+,s,t)ds.\label{eq1.22}
\end{equation}

\begin{lemma}\label{l1.1}
For the kernels $K_n(z,x,t)$, the following estimations hold:
\begin{equation}
|K_n(z,x,t)|<d(z(t-x))\frac{\sigma_+^{n-1}(x)}{(n-1)!},\label{eq1.23}
\end{equation}
where
\begin{equation}
d(z)=e^{|\beta|}\cosh\frac{\alpha\sqrt3}2\,(z=\alpha+i\beta;\alpha,\beta\in\mathbb{R});\quad\sigma_+(x)=\int\limits_x^\infty\left|\frac{m(t)}{m_+}-1\right|dt.\label{eq1.24}
\end{equation}
\end{lemma}

P r o o f. Equalities
$$s_2(iz(x-t))=\frac13\left\{e^{iz\zeta_0(x-t)}+\frac1{\zeta_1^2}e^{iz\zeta_1(x-t)}+\frac1{\zeta_2^2}e^{iz\zeta_2(x-t)}\right\};$$
$$iz\zeta_0=i\alpha-\beta;\, iz\zeta_1=\frac12[(\beta+\alpha\sqrt3)+i(-\alpha-\beta\sqrt3)];\, iz\zeta_2=\frac12[(\beta-\alpha\sqrt3)+i(-\alpha+\beta\sqrt3)]$$
imply that, for the kernel $K_1(z,x,t)$ \eqref{eq1.20},
$$|K_1(z,x,t)|<d(z(t-x))\quad(t>x),$$
where $d(z)$ is given by \eqref{eq1.24}, which coincides with \eqref{eq1.23} when $n=1$. Using induction by $n$, equality \eqref{eq1.22}, and estimations \eqref{eq1.23}, we obtain
$$|K_{n+1}(z,x,t)|\leq\int\limits_x^td(z(s-x))\left|\frac{m(s)}{m_+}-1\right|d(z(s-x))\frac{\sigma_+^{n-1}(s)}{(n-1)!}ds.$$
And since
$$d(z(s-x))\cdot d(z(s-t))=4e^{|\beta|(t-x)}\cosh\frac{\alpha\sqrt3}2(s-x)\cosh\frac{\alpha\sqrt3}2(t-s)$$
$$=e^{|\beta|(t-x)}\frac42\left\{\cosh\frac{\alpha\sqrt3}2(t-x)+\cosh\frac{\alpha\sqrt3}2(t-x-2s)\right\}<d(z(t-x)),$$
then
$$|K_{n+1}(z,x,t)|\leq d(z(t-x))\frac{\sigma_+^n(x)}{n!},$$
this proves \eqref{eq1.23} for $n+1$. $\blacksquare$

Relation \eqref{eq1.19} yields that
\begin{equation}
v_k(\lambda,x)=e^{iz_+\zeta_kx}+\int\limits_x^\infty N(z_+,x,t)\left(\frac{m(t)}{m_+}-1\right)e^{iz_+\zeta_kt}dt\label{eq1.25}
\end{equation}
where
\begin{equation}
N(z,x,t)=\sum\limits_{n=1}^\infty(-iz)^nK_n(z,x,t).\label{eq1.26}
\end{equation}
Series \eqref{eq1.26} is majorized by the converging (uniformly by $x$) series, due to \eqref{eq1.23},
\begin{equation}
|N(z,x,t)|<d(z(t-x))\exp\{|z|\sigma_+(x)\}.\label{eq1.27}
\end{equation}
Proceed to the boundary condition (a) \eqref{eq1.16}. Relations \eqref{eq1.25} -- \eqref{eq1.27} imply
$$|v_k(\lambda,x)|\leq e^{|z|x}\left(1+\exp\{|z|\sigma_+\}\int\limits_x^\infty e^{2|z|(t-x)}\left|\frac{m(t)}{m_+}-1\right|dt\right)$$
where $\sigma_+=\sigma_+(0)$. And since ${\displaystyle e^{2|z|(t-x)}<e^{2|z|t}}$ ($t>x>0$), then for $2|z|<a$, the last integral converges, due to (ii) \eqref{eq1.12}. So,
$$|v_k(\lambda,x)|<e^{|z|x}(1+\exp\{|z|\sigma_+\}\cdot A_+(z));\quad A_+(z)\stackrel{\rm def}{=}\int\limits_0^\infty e^{2|z|t}\left|\frac{m(t)}{m_+}-1\right|dt.$$
Finally, equation \eqref{eq1.19} implies that
$$|v_k(\lambda,x)-e^{iz_+\zeta_kx}|<|z|(1+\exp|z|\sigma_+)A_+(z)\int\limits_x^\infty e^{2|z|t}\left|\frac{m(t)}{m_+}-1\right|dt.$$
The last expression after the integral sign $\displaystyle{e^{2|z|t}\left|\frac{m(t)}{m_+}-1\right|\in L^1(\mathbb{R}_+)}$, for $2|z|<a$. And thus, for every such $z$ ($2|z|<a$) condition (a) \eqref{eq1.16} holds.

Define the disk
\begin{equation}
\mathbb{D}_{a/2}=\{z\in\mathbb{C}:|z|<a/2\},\label{eq1.28}
\end{equation}
and, in terms of $\lambda=z_+/n_+$,
\begin{equation}
\mathbb{D}_{a/2}^+=\left\{\lambda\in\mathbb{C}:|\lambda|<\frac a{2n_+}\right\}.\label{eq1.29}
\end{equation}

\begin{theorem}\label{t1.1}
For all $\lambda\in\mathbb{D}_{a/2}^+$ \eqref{eq1.29}, the Jost solutions $\{v_k(\lambda,x)\}_0^2$ to the boundary value problem \eqref{eq1.11}, (a) \eqref{eq1.16} exist and are given by \eqref{eq1.25}.
\end{theorem}

Analogously to \eqref{eq1.19}, the Jost solutions $\{u_k(\lambda,x)\}_0^2$ to the boundary value problem \eqref{eq1.11}, (b) \eqref{eq1.16} are matched with the integral equation
\begin{equation}
u_k(\lambda,x)=e^{iz_-\zeta_kx}+iz_-\int\limits_{-\infty}^xs_2(iz_-(x-t))\left(\frac{m(t)}{m_-}-1\right)dt\quad(0\leq k\leq2)\label{eq1.30}
\end{equation}
where $z_-=\lambda n_-$. Equation \eqref{eq1.30}, in terms of the Volterra operator
$$(\widehat{K}_zf)(x)=\int\limits_{-\infty}^x\widehat{K}_1(z,x,t)\left(\frac{m(t)}{m_-}-1\right)f(t)dt,$$
where
$$\widehat{K}_1(z,x,t)=s_2(iz(x-t)),$$
becomes
\begin{equation}
(I-iz_-\widehat{K}_{z_-})u_k(\lambda,x)=e^{iz_-\zeta_kx},\label{eq1.31}
\end{equation}
therefore
$$u_k(\lambda,x)=\sum(iz_-)^n\widehat{K}_{z_-}^ne^{iz_-\zeta_kx}.$$
The operators $\widehat{K}_z^n$ are also Volterra,
$$(\widehat{K}_z^nf)(x)=\int\limits_{-\infty}^x\widehat{K}_n(z,x,t)\left(\frac{m(t)}{m_-}-1\right)f(t)dt,$$
and the recurrent relations hold for the kernels $\widehat{K}_n(z,x,t)$:
$$\widehat{K}_{n+1}(z,x,t)=\int\limits_{-\infty}^x\widehat{K}_n(z,x,s)\left(\frac{m(s)}{m_-}-1\right)\widehat{K}_1(z,s,t)ds.$$
Analogously to Lemma \ref{l1.1}, the following statement is true.

\begin{lemma}\label{l1.2}
For the kernels $\widehat{K}_n(z,x,t)$, the following estimates hold:
\begin{equation}
|\widehat{K}_n(z,x,t)|<d(z(x-t))\frac{\sigma_-^{n-1}(x)}{(n-1)!}\label{eq1.32}
\end{equation}
where
\begin{equation}
d(z)=e^{|\beta|}\cosh\frac{\alpha\sqrt3}2\quad(z=\alpha+i\beta,\alpha,\beta\in\mathbb{R});\quad\sigma_-(x)=\int\limits_{-\infty}^x\left|\frac{m(t)}{m_-}-1\right|dt.\label{eq1.33}
\end{equation}
\end{lemma}

Equation \eqref{eq1.31} implies that
\begin{equation}
u_k(\lambda,x)=e^{iz_-\zeta_kx}+\int\limits_{-\infty}^x\widehat{N}(z_-,x,t)\left(\frac{m(t)}{m_-}-1\right)e^{iz_-\zeta_kt}dt\label{eq1.34}
\end{equation}
where
\begin{equation}
\widehat{N}(z,x,t)\stackrel{\rm def}{=}\sum(iz)^nK_n(z,x,t).\label{eq1.35}
\end{equation}
Due to \eqref{eq1.32}, this series is majorized by the convergent (uniformly by $x$) series and
$$|\widehat{N}(\lambda,x,t)|<dz(t-x)\cdot\exp\{|z|\sigma_t(x)\}.$$
Define the disk $\mathbb{D}_{a/2}^-$ obtained from $\mathbb{D}_{a/2}$ \eqref{eq1.22} as a result of the substitution $z=\lambda n_-$,
\begin{equation}
\mathbb{D}_{a/2}^-\stackrel{\rm def}{=}\left\{\lambda\in\mathbb{C}:|\lambda|<\frac a{2n_-}\right\}.\label{eq1.36}
\end{equation}

\begin{theorem}\label{t1.2}
For all $\lambda\in\mathbb{D}_{a/2}^-$ \eqref{eq1.36}, Jost solutions $\{u_k(\lambda,x)\}_0^2$ to the boundary value problem \eqref{eq1.11}, (b) \eqref{eq1.16} exist and are given by \eqref{eq1.34}.
\end{theorem}

Analogously to the Wronskian of two functions $u(x)$, $v(x)$
$$W(u,v)\stackrel{\rm def}{=}\det\left|
\begin{array}{ccc}
u(x)&v(x)\\
u'(x)&v'(x)
\end{array}\right|$$
linked with a second-order equation, we define the Wronskian for three functions $u(x)$, $v(x)$, $w(x)$ (the {\bf third-order Wronskian}) by the formula
\begin{equation}
W(u,v,w)\stackrel{\rm def}{=}\det\left[
\begin{array}{ccc}
u(x)&v(x)&w(x)\\
u'(x)&v'(x)&w'(x)\\
u''(x)&v''(x)&w''(x)
\end{array}\right].\label{eq1.37}
\end{equation}
This Wronskian corresponds to a third-order equation.

Jost solutions $\{v_k(\lambda,x)\}_0^2$ (as well as $\{u_k(\lambda,x)\}_0^2$) are linearly independent. To prove this, calculate the Wronskian $W(v_0,v_1,v_2)$. First of all, notice that ${\displaystyle\frac d{dx}W(v_0,v_1,v_2)=0}$, in view of equation \eqref{eq1.11}. Using (a) \eqref{eq1.16}, we have
$$W(v_0,v_1,v_2)=\det\left[
\begin{array}{ccc}
e^{iz_+\zeta_0x}&e^{iz_+\zeta_1x}&e^{iz_+\zeta_2x}\\
(iz_+\zeta_0)e^{iz_+\zeta_0x}&(iz_+\zeta_1)e^{iz_+\zeta_1x}&(iz_+\zeta_2)e^{iz_+\zeta_2x}\\
(iz_+\zeta_0)^2e^{iz_+\zeta_0x}&(iz_+\zeta_1)^2e^{iz_+\zeta_1x}&(iz_+\zeta_2)^2e^{iz_+\zeta_2x}
\end{array}\right]$$
$$=m_+\lambda^3\det\left[
\begin{array}{ccc}
1&1&1\\
1&\zeta_1&\zeta_2\\
1&\zeta_2&\zeta_1
\end{array}\right]=3\sqrt3m_+\lambda^3.$$
So,
\begin{equation}
W(v_0,v_1,v_2)=3\sqrt3m_+\lambda^3.\label{eq1.38}
\end{equation}
Hence follows the statement.

\begin{lemma}\label{l1.3}
For all $\lambda\not=0$, functions $\{v_k(\lambda,x)\}$ are linearly independent.
\end{lemma}

\begin{remark}\label{r1.3}
The Wronskian of Jost solutions $\{u_k(\lambda,x)\}_0^2$ equals
\begin{equation}
W(u_0,u_1,u_2)=3\sqrt3m_+\lambda^3.\label{eq1.39}
\end{equation}
\end{remark}
\vspace{5mm}

{\bf 1.3} Turn to the analytic (by $\lambda$) properties of the functions $\{v_k(\lambda,x)\}_0^2$. Define the functions
\begin{equation}
\psi_k(\lambda,x)=v_k(\lambda,x)e^{-iz_+\zeta_kx}\quad(z_+=\lambda n_+,0\leq k\leq2),\label{eq1.40}
\end{equation}
which, due to \eqref{eq1.19} are solutions to the equations
\begin{equation}
\psi_k(\lambda,x)=1-iz_+\int\limits_x^\infty e^{iz_+\zeta_k(t-x)}s_2(iz_+(x-t))\left(\frac{m(t)}{m_+}-1\right)\psi_k(\lambda,t)dt\quad(0\leq k\leq2).\label{eq1.41}
\end{equation}
Equation \eqref{eq1.19} implies two more equalities:
\begin{equation}
\begin{array}{llll}
(i)\,v'_k(\lambda,x)e^{-iz_+\zeta_kx}\\
=iz_+\left[{\displaystyle\zeta_k-iz_+\int\limits_x^\infty e^{iz_+\zeta_k(t-x)}s_1(iz_+(x-t))\left(\frac{m(t)}{m_+}-1\right)\psi_k(\lambda,t)}\right]\,(0\leq k\leq2);\\
(ii)\,v''_k(\lambda,x)e^{-iz_+\zeta_kx}\\
=(iz_+)^2\left[\zeta_k^2-iz_+{\displaystyle\int\limits_x^\infty e^{iz_+\zeta_k(t-x)}s_0(iz_+(x-t))\left(\frac{m(t)}{m_+}-1\right)\psi_k(\lambda,t)dt}\right]\,(0\leq k\leq2).
\end{array}\label{eq1.42}
\end{equation}

By $\{S_p(i)\}_0^5$, we define the sectors obtained from $\{S_p\}_0^5$ \eqref{eq1.10} after the rotation through $\pi/2$,
\begin{equation}
S_p(i)=\{i\lambda:\lambda\in\ S_p\}\quad(0\leq p\leq5).\label{eq1.43}
\end{equation}

\begin{picture}(200,200)
\put(0,100){\vector(1,0){200}}
\put(100,0){\vector(0,1){200}}
\put(0,133){\vector(3,-1){200}}
\put(200,133){\vector(-3,-1){200}}
\put(105,185){$il_{\zeta_0}$}
\put(120,150){$\Omega_0$}
\qbezier(153,120)(100,140)(49,120)
\qbezier(100,23)(155,35)(168,125)
\qbezier(45,120)(30,70)(100,25)
\put(140,105){$S_4(i)$}
\put(70,110){$S_0(i)$}
\put(70,70){$S_2(i)$}
\put(50,100){$S_1(i)$}
\put(30,40){$\Omega_2$}
\put(0,50){$il_{\zeta_1}$}
\put(115,70){$S_3(i)$}
\put(105,110){$S_5(i)$}
\put(180,80){$\Omega_1$}
\put(185,105){$l_{\zeta_0}$}
\put(180,60){$il_{\zeta_2}$}
\end{picture}

\hspace{20mm} Fig. 1

Define three sectors:
$$\Omega_0\stackrel{\rm def}{=}S_5(i)\cup S_0(i)\cup (il_{\zeta_0});$$
\begin{equation}
\Omega_1\stackrel{\rm def}{=}S_2(i)\cup S_4(i)\cup(il_{\zeta_2});\label{eq1.44}
\end{equation}
$$\Omega_2\stackrel{\rm def}{=}S_1(i)\cup S_2(i)\cup(il_{\zeta_1}).$$
By $\{\Omega_k^-\}_0^2$, we denote sectors that are centrally symmetric to the sectors $\{\Omega_k\}_0^2$ \eqref{eq1.44},
\begin{equation}
\Omega_k^-\stackrel{\rm def}{=}\{\lambda:-\lambda\in\Omega_k\}.\label{eq1.45}
\end{equation}
At the relations \eqref{eq1.42}, \eqref{eq1.43}, consider the case of $k=0$, cases of $k=1$ and $k=2$ follow from $k=0$ after the substitutions $\lambda\rightarrow\lambda\zeta_1$, $\lambda\rightarrow\lambda\zeta_2$.

First, we study behavior of the kernel $e^{iz_+\zeta_0(t-x)}s_2(iz_+(x-t))$.

\begin{lemma}\label{l1.4}
For all $\tau>0$, the following equalities are true:
\begin{equation}
(i)\,e^{iz_+\zeta_k\tau}s_2(-iz_+\tau)=\frac{\zeta_k}3+r_k(\lambda)\quad(\lambda\in\Omega_k,|\lambda|\gg1,0\leq k\leq2),\label{eq1.46}
\end{equation}
besides, $r_k(\lambda)\lambda^p\rightarrow0$ as $\lambda\rightarrow\infty$ ($\lambda\in\Omega_k$ and for all $p\in\mathbb{Z}_+$);
\begin{equation}
(ii)\,e^{-iz_-\zeta_k\tau}s_2(-iz_-\tau)=\frac{\zeta_k}3+\widehat{r}_k(\lambda)\quad(\lambda\in\Omega_k^-,|\lambda|\gg1,0\leq k\leq2),\label{eq1.47}
\end{equation}
besides, $\widehat{r}_k(\lambda)\lambda^p\rightarrow0$ when $\lambda\rightarrow\infty$ ($\lambda\in\Omega_k^-$, $\forall p\in\mathbb{Z}_+$).
\end{lemma}

P r o o f. Consider the case of $k=0$. Using
\begin{equation}
e^{iz_+\zeta_0\tau}s_2(-iz_+\tau)=\frac13\left[1+\frac1{\zeta_1^2}e^{i(\zeta_0-\zeta_1)z_+\tau}+\frac1{\zeta_2^2}e^{i(\zeta_0-\zeta_2)z_+\tau}\right]=\frac13+r_0(\lambda),\label{eq1.48}
\end{equation}
also $z_+=n_+\lambda$ ($\lambda=\alpha+i\beta$), and
$$i(\zeta_0-\zeta_1)z_+=\frac{\sqrt3n_+}2[(\alpha-\beta\sqrt3)+i(\alpha\sqrt3+\beta)];$$
$$i(\zeta_0-\zeta_2)z_+=\frac{n_+\sqrt3}2[(-\alpha-\beta\sqrt3)+i(\alpha\sqrt3-\beta)],$$
we obtain that for $\tau>0$ in the domain $\{\lambda=\alpha+i\beta:\alpha<\beta\sqrt3,-\alpha<\beta\sqrt3\}$, which exactly coincides with $\Omega_0$ \eqref{eq1.44}, the following estimates hold:
\begin{equation}
\begin{array}{ccc}
{\displaystyle|e^{i(\zeta_0-\zeta_1)z_+\tau}|<\exp\frac{\sqrt3}2n_+\tau(\alpha-\beta\sqrt3);}\\
{\displaystyle|e^{i(\zeta_0-\zeta_2)z_+\tau}|<\exp\frac{\sqrt3}2n_+\tau(-\alpha-\beta\sqrt3).}
\end{array}\label{eq1.49}
\end{equation}
Each of these exponents is less than one, for $\lambda\in\Omega_0$, and vanishes while exponentially decreasing  as $\lambda\rightarrow\infty$ ($\lambda\in\Omega_0$), and thus $r_0(\lambda)\lambda^p\rightarrow0$ ($p\in\mathbb{Z}_+$) when $\lambda\rightarrow\infty$. Other equalities in \eqref{eq1.48} are proved analogously. Equalities \eqref{eq1.48} follow from \eqref{eq1.47}. $\blacksquare$

\begin{remark}\label{r1.2}
For the kernels $e^{iz_+\zeta_k\tau}s_1(-z_+\tau)$ and $e^{iz_+\zeta_k\tau}s_0(iz_+\tau)$, analogues of the equalities \eqref{eq1.48} are true in the corresponding sectors $\{\Omega_k\}_0^2$, and for the kernels $e^{-iz_-\zeta_k\tau}s_1(-z_-\tau)$ and $e^{-iz\zeta_k\tau}s_0(-z_-\tau)$, analogues of the relations \eqref{eq1.49} hold in the sectors $\{\Omega_k^-\}_0^2$ \eqref{eq1.45}.
\end{remark}

Consider equation \eqref{eq1.42} when $k=0$,
\begin{equation}
\psi_0(\lambda,x)=1-iz_+\int\limits_x^\infty\widetilde{K}_1(z_+,x,t)\left(\frac{m(t)}{m_+}-1\right)\psi_0(\lambda,t)dt\label{eq1.50}
\end{equation}
where
\begin{equation}
\widetilde{K}_1(z,x,t)=e^{iz\zeta_(t-x)}s_2(z_+(x-t)).\label{eq1.51}
\end{equation}
In terms of the Volterra operator
$$(\widetilde{T}_zf)(x)=\int\limits_x^\infty\widetilde{K}_1(z,x,t)\left(\frac{m(t)}{m_+}-1\right)f(t)dt,$$
where $\widetilde{K}_1(z,x,t)$ is given by \eqref{eq1.51}, equation \eqref{eq1.51} becomes
$$(I+iz_+\widetilde{K}_{z_+})\psi_0(\lambda,x)=1,$$
and thus
\begin{equation}
\psi_0(\lambda,x)=\sum\limits_0^\infty(+iz_+)^nK_{z_+}^n1.\label{eq1.52}
\end{equation}
Operators $\widetilde{K}_z^n$ are also Volterra,
$$(\widetilde{K}_z^nf)(x)=\int\limits_x^\infty\widetilde{K}_n(z,x,t)\left(\frac{m(t)}{m_+}-1\right)f(t)dt,$$
and their kernels $\widetilde{K}_n(z,x,t)$ satisfy the recurrent relation
\begin{equation}
\widetilde{K}_{n+1}(z,x,t)=\int\limits_x^t\widetilde{K}_1(z,x,s)\left(\frac{m(s)}{m_+}-1\right)\widetilde{K}_n(z,s,t)ds.\label{eq1.53}
\end{equation}
For the kernel $\widetilde{K}_1(z,x,t)$ \eqref{eq1.52}, due to \eqref{eq1.48}, \eqref{eq1.49}, we have the estimation
\begin{equation}
|\widetilde{K}_1(z,x,t)|<\widetilde{d}_1(z(t-x));\quad \widetilde{d}_1(z(t-x))\stackrel{\rm def}{=}\frac13(1+2e^{-\beta\sqrt3(t-x)}\cosh\alpha(t-s)).\label{eq1.54}
\end{equation}

\begin{lemma}\label{l1.5}
For all $\lambda\in\Omega_0$, the following inequalities for the kernels $\widetilde{K}_n(z,x,t)$ are true:
\begin{equation}
|\widetilde{K}_n(z,x,t)|<\widetilde{d}_n(z,x,t)\frac{\sigma_+^{n-1}(x)}{(n-1)!},\label{eq1.55}
\end{equation}
besides,
\begin{equation}
\begin{array}{ccc}
{\displaystyle\widetilde{d}_n(z,x,t)\stackrel{\rm def}{=}\frac1{3^n}\{a_n+2^{n-1}\cdot e^{-\beta\sqrt3(t-x)}\cosh\alpha(t-x)\}}\quad(z=\alpha+i\beta\in\Omega_0,\,\alpha,\beta\in\mathbb{R});\\
{\displaystyle a_{n+1}<3a_n+2^n\quad(a_1=1,n\geq2)\quad\sigma_+(x)\stackrel{\rm def}{=}\int\limits_x^\infty\left|\frac{m(t)}{m_+}-1\right|dt.}
\end{array}\label{eq1.56}
\end{equation}
\end{lemma}

P r o o f. Relations \eqref{eq1.53} and \eqref{eq1.54}, \eqref{eq1.55} imply
\begin{equation}
|\widetilde{K}_{n+1}(z,x,t)|<\int\limits_x^t|\widetilde{d}_1(z,x,s)\widetilde{d}_n(z,s,t)|\left|\frac{m(s)}{m_+}-1\right|\frac{\sigma_+^{n-1}(s)}{(n-1)!}ds.\label{eq1.57}
\end{equation}
Taking into account \eqref{eq1.55} and \eqref{eq1.56} for $n$, we obtain that
$$|\widetilde{d}_1(z,x,s)\cdot\widetilde{d}_n(z,s,t)|<\frac1{3^{n+1}}(1+2e^{-i\beta\sqrt3(s-x)}\cosh\alpha(s-x))$$
$$\times(a_n-2^ne^{-\beta\sqrt3(t-s)}\cosh\alpha(t-s))$$
$$=\frac1{3^{n+1}}(a_n+a_n2e^{-\beta\sqrt3(s-x)}\cosh\alpha(s-x)+2e^{-\beta\sqrt3(t-s)}\cosh\alpha(t-s))$$
$$+4e^{-\sqrt3(t-x)}\cosh\alpha(t-s)\cosh\alpha(s-x).$$
Notice that
$$\cosh\alpha(s-x)\cosh\alpha(t-s)=\frac2z(\cosh\alpha(t-x)+\cosh\alpha(x+t-2s))<\cosh\alpha(t-s).$$
Therefore,
$$|\widetilde{d}_1(z,x,s)\widetilde{d}_n(z,s,t)|<\frac1{3^{n+1}}\{a_n+2^{n-1}\cdot e^{-\beta\sqrt3(t-s)}\cosh\alpha(t-s)+a_n2\cdot e^{-\beta(s-x)}$$
$$\times\cosh\alpha(s-x)+2^ne^{-\beta\sqrt3(t-x)}\cosh\alpha(t-x)\}.$$
Function $f(s)=2e^{-\beta\sqrt3(t-s)}\cosh\alpha(t-s)$ on the interval $s\in(x,t)$ monotonically increases from $f(x)=e^{-\beta\sqrt3(t-x)}\cosh\alpha(t-x)$ to 2, function $f_1(s)=e^{-\beta\sqrt3(s-x)}$ $\times\cosh\alpha(s-x)$ on the same interval $s\in(x,t)$ decreases from 2 to $f_1(t)=2e^{-\beta\sqrt3(t-x)}\cosh\alpha(t-s).$ As a result, we have
$$|d_1(z,x,s)\widetilde{d}_n(z,x,t)|<\frac1{3^{n+1}}\{a_n+2a_n+2^n+2^ne^{-\beta\sqrt3(t-x)}\cosh\alpha(t-x)\}$$
$$<\frac1{3^{n+1}}\{a_{n+1}+2^ne^{-\beta\sqrt3(t-x)}\cosh\alpha(t-x)\}.$$
Using this inequality from \eqref{eq1.57}, we obtain that
$$|K_{n+1}(z,x,t)|<d_{n+1}(z,x,t)\frac{\sigma_+^n(x)}{n!},$$
which proves \eqref{eq1.56} for $n+1$. Therefore, due to the principle of mathematical induction, \eqref{eq1.56} is proved  for all $n\in\mathbb{N}$. $\blacksquare$

\begin{remark}\label{r1.3}
The inequality $a_{n+1}<3a_n+2^n$ ($n\in\mathbb{N}$) implies that
$$a_{n+1}<3^{n+1}-2^{n+1}\quad(\forall n\in\mathbb{Z}_+).$$
Hence, for $d_n(z,x,t)$ \eqref{eq1.57}, we obtain the (uniform by $n$) estimate
\begin{equation}
d_n(z,x,t)<1+e^{-\beta\sqrt3(t-x)}\cosh\alpha(t-x).\label{eq1.58}
\end{equation}
\end{remark}

Using \eqref{eq1.53}, we obtain
\begin{equation}
\psi_0(\lambda,x)=1+iz_+\int\limits_x^\infty\widetilde{N}(z,x,t)\left(\frac{m(t)}{m_+}-1\right)1dt,\label{eq1.59}
\end{equation}
where
\begin{equation}
\widetilde{N}(z,x,t)\stackrel{\rm def}{=}\sum\limits_1^\infty\widetilde{K}_n(z,x,t).\label{eq1.60}
\end{equation}
In view of \eqref{eq1.55} and \eqref{eq1.58}, series \eqref{eq1.60} is majorized by the converging series,
$$|\widetilde{N}(z,x,t)|<(1+e^{-\beta\sqrt3(t-x)}\cosh\alpha(x-t))\exp|z|\sigma_+(x)\quad(\forall\lambda\in\Omega_0).$$
Hence we obtain the estimate for $\psi_0(\lambda,x)$ \eqref{eq1.59},
\begin{equation}
|\psi_0(\lambda,x)|<1+e^{|z_+|\sigma_+(x)}<2e^{|z_+|\sigma_+(x)}<2e^{|z_+|\sigma_+(0)}\quad(\forall\lambda\in\Omega_0),\label{eq1.61}
\end{equation}
and thus, for all $\lambda\in\Omega_0$, function $\psi_0(\lambda,x)$ is bounded for all $x\in\mathbb{R}_+$.

\begin{theorem}\label{t1.4}
For all $\lambda\in\Omega_0$, functions $v_0(\lambda,x)$, $v'_0(\lambda,x)$, $v''_0(\lambda,x)$ belong to the space $L^2(\mathbb{R}_+)$.
\end{theorem}

P r o o f. Since
\begin{equation}
v_0(\lambda,x)=e^{iz_+x}\psi_0(\lambda,x),\label{eq1.62}
\end{equation}
then the boundedness of $\psi_0(\lambda,x)$ ($\forall\lambda\in\Omega_0$) and inclusion $e^{iz_+x}\in L^2(\mathbb{R}_+)$ ($\forall\lambda\in\Omega_0$) imply that $v_0(\lambda,x)\in L^2(\mathbb{R}_+)$.

Equation (i) \eqref{eq1.42} yields that
\begin{equation}
v'_0(\lambda,x)=iz_+e^{iz_+x}\psi_0^1(\lambda,x)\label{eq1.63}
\end{equation}
where
$$\psi_0^1(\lambda,x)=1-iz_+\int\limits_x^\infty e^{iz_+(t-x)}s_1(iz_+(x-t))\left(\frac{m(t)}{m_+}-1\right)\psi_0(\lambda,t).$$
Due to \eqref{eq1.61}, function $\psi_0^1(\lambda,x)$  is bounded,
$$|\psi_0(\lambda,x)|<3e^{|z_+|\sigma_+(x)}<3e^{|z_+|\sigma_+(0)}\quad(\forall\lambda\in\Omega_0),$$
which gives us that $v'(\lambda,x)\in L^2(\mathbb{R}_+)$ ($\forall\lambda\in\Omega_0$).

Using (ii) \eqref{eq1.42}, we have
\begin{equation}
\psi''_0(\lambda,x)=(iz_+)^2e^{iz_+x}\psi_0^2(\lambda,x),\label{eq1.64}
\end{equation}
where
$$\psi_0^2(\lambda,x)=1-iz_+\int\limits_x^\infty e^{iz_+(t-x)}s_0(iz_+(x-t))\left(\frac{m(t)}{m_+}-1\right)\psi_0(\lambda,t)dt,$$
besides, function $\psi_0^2(\lambda,x)$ is also bounded,
$$|\psi_0^2(\lambda,x)|<3e^{|z_+|\sigma_+(x)}<3e^{|z_+|\sigma_+(0)}\quad(\forall\lambda\in\Omega_0).$$
Therefore, \eqref{eq1.64} implies that $v''_0(\lambda,x)\in L^2(\mathbb{R}_+)$ ($\forall\lambda\in\Omega_0$). $\blacksquare$

\begin{remark}\label{r1.4}
For the functions $v_1(\lambda,x)$ and $v_2(\lambda,x)$, statement of Theorem \ref{t1.4} holds for all $\lambda\in\Omega_1$ and for all $\lambda\in\Omega_2$ correspondingly. Moreover, analogues of the formulas \eqref{eq1.62} -- \eqref{eq1.64} are true.
\end{remark}

\begin{theorem}\label{t1.5}
Solution $\psi_0(\lambda,x)$ to equation \eqref{eq1.41} ($k=0$) is a holomorphic function in the sector $\Omega_0$, and the following relations hold:\\
(i) $\lim\limits_{x\rightarrow\infty}\psi_0(\lambda,x)=1$ ($\lambda\in\Omega_0$);\\
(ii) $\lim\limits_{\lambda\rightarrow\infty}\psi_0(\lambda,x)=1$ ($\lambda\in\Omega_0,x\in\mathbb{R}_+$);\\
(iii) $\displaystyle{\lim\limits_{\lambda\rightarrow\infty}\frac{\psi_0(\lambda,x)-1}{(iz_+)^3}=M_+(x)}$; $\displaystyle{M_+(x)\stackrel{\rm def}{=}\int\limits_x^\infty\frac{(x-t)^2}2\left(\frac{m(t)}{m_+}-1\right)dt}$
\begin{equation}
(\lambda\in\Omega_0,x\in\mathbb{R}_+);\label{eq1.65}
\end{equation}
(iv) ${\displaystyle\frac{d^3}{dx^3}M_+(x)=(1-x)\left(\frac{m(x)}{m_t}-1\right)}$ ($\forall x\in\mathbb{R}_+$).
\end{theorem}
\vspace{5mm}

{\bf 1.4} Functions
\begin{equation}
\varphi_k(\lambda,x)=u_k(\lambda,x)e^{-iz_-\zeta_kx}\quad(0\leq k\leq2),\label{eq1.66}
\end{equation}
due to \eqref{eq1.30}, are solutions to the integral equations
\begin{equation}
\varphi_k(\lambda,x)=1+(iz_-)\int\limits_{-\infty}^xe^{iz_-\zeta_k(t-x)}s_2(iz_-(x-t))\left(\frac{m(t)}{m_-}-1\right)\varphi_k(\lambda,t)dt\quad(0\leq k\leq2)\label{eq1.67}
\end{equation}
and analogously to \eqref{eq1.42}
\begin{equation}
\begin{array}{llll}
(i)u'_k(\lambda,x)e^{-iz_-\zeta_kx}\\
=\displaystyle{iz_-\zeta_k+(iz_-)^2\int\limits_{-\infty}^xe^{i\zeta_kz_-(t-x)}s_1(iz_-(x-t))\left(\frac{m(t)}{m_-}-1\right)\varphi_k(\lambda,t)dt;}\\
(ii)u''_k(\lambda,x)e^{-iz_-\zeta_kx}\\
\displaystyle{=(iz_-\zeta_k)^2+(iz_-)^3\int\limits_{-\infty}^xe^{iz_-\zeta_k(t-x)}s_0(iz_-(x-t))\left(\frac{m(t)}{m_-}-1\right)\varphi_k(\lambda,t)dt.}
\end{array}\label{eq1.68}
\end{equation}
($0\leq k\leq2$). The following statement is an analogue of Theorem \ref{t1.4}.

\begin{theorem}\label{t1.6}
For all $\lambda\in\Omega_0^-$, the functions $u_0(\lambda,x)$, $u'_0(\lambda,x)$, and $u''_0(\lambda,x)$ belong to the space $L^2(\mathbb{R}_-)$.
\end{theorem}

P r o o f. Analogously to \eqref{eq1.61}, for $\varphi_0(\lambda,x)$,
\begin{equation}
|\varphi_0(\lambda,x)|<2e^{|z_-|\sigma_-(x)}<2e^{|z_-|\sigma_-(0)}\label{eq1.69}
\end{equation}
when $\lambda\in\Omega_0^-$. Using equality
\begin{equation}
u_0(\lambda,x)=e^{iz_-x}\varphi_0(\lambda,x),\label{eq1.70}
\end{equation}
and $e^{iz_-x}\in L^2(\mathbb{R}_-)$ for all $\lambda\in\Omega_0^-$, and also \eqref{eq1.69}, we obtain that $u_0(\lambda,x)\in L^2(\mathbb{R}_-)$ for all $\lambda\in\Omega_0^-$. Equation (i) \eqref{eq1.68} implies
\begin{equation}
u'_0(\lambda,x)=iz_-e^{iz_-x}\varphi_0^1(\lambda,x),\label{eq1.71}
\end{equation}
where
$$\varphi_0^1(\lambda,x)=1+iz_-\int\limits_{-\infty}^xe^{iz_-(t-x)}s_1(iz_-(x-t))\left(\frac{m(t)}{m_-}-1\right)\varphi_0(\lambda,t)dt.$$
And since
$$|\varphi_0^1(\lambda,x)|<3e^{|z_-|\sigma_-(x)}<3e^{|z_-|\sigma_-(0)}\quad(\forall\lambda\in\Omega_0^-),$$
then \eqref{eq1.71} and $e^{iz_-x}\in L^2(\mathbb{R}_-)$ ($\forall\lambda\in\Omega_0^-$) imply that $u'_0\in L^2(\mathbb{R}_-)$ ($\forall\lambda\in\Omega_0^8$). Using (ii) \eqref{eq1.68}, we have
\begin{equation}
u''_0(\lambda,x)=(iz_-)^2e^{iz_-x}\varphi_0^2(\lambda,x)\label{eq1.72}
\end{equation}
where
$$\varphi_0^2(\lambda,x)=1+iz_-\int\limits_{-\infty}^xe^{iz_-(t-x)}s_0(iz_-(x-t))\left(\frac{m(t)}{m_-}-1\right)\varphi_0(\lambda,t)dt$$
and, taking into account that
$$|\varphi_0^2(\lambda,x)|<3e^{|z_-|\sigma_-(x)}<3e^{|z_-|\sigma_0(0)},$$
equation \eqref{eq1.72} yields that $u''_0(\lambda,x)\in L^2(\mathbb{R}_-)$ for all $\lambda\in\Omega_0^-$. $\blacksquare$

\begin{remark}\label{r1.5}
For the functions $u_1(\lambda,x)$ and $u_2(\lambda,x)$, Theorem \ref{t1.6} and analogues of the formulas \eqref{eq1.70} -- \eqref{eq1.72} also hold in the sectors $\Omega_1^-$ and $\Omega_2^+$.
\end{remark}

\begin{theorem}\label{t1.7}
Solution $\varphi_0(\lambda,x)$ to equation \eqref{eq1.67} $k=0$ is an analytic function from $\lambda$ in the sector $\Omega_0^-$ and\\
(i) $\lim\limits_{x\rightarrow-\infty}\varphi_0(\lambda,x)=0$ ($\lambda\in\Omega_0^-$);\\
(ii) $\lim\limits_{\lambda\rightarrow\infty}\varphi_0(\lambda,x)=1$ ($\lambda\in\Omega_0^-$, $x\in\mathbb{R}_-$);\\
(iii) ${\displaystyle\lim\limits_{\lambda\rightarrow\infty}\frac{\varphi_0(\lambda,x)-1}{(iz_-)^3}=M_-(x)}$; $\displaystyle{M_-(x)\stackrel{\rm def}{=}\int\limits_{-\infty}^x\frac{(x-t)}2\left(\frac{m(t)}{m_-}-1\right)dt}$,
\begin{equation}
(\lambda\in\Omega_0^-,x\in\mathbb{R});
\label{eq1.73}
\end{equation}
(iv) ${\displaystyle\frac{d^3}{dx^3}M_-(x)=(1-x)\left(\frac{m(x)}{m_-}-1\right)}$ for all $x\in\mathbb{R}_-$.
\end{theorem}

\begin{remark}\label{r1.6}
For the functions $\psi_1(\lambda,x)$ and $\psi_2(\lambda,x)$ \eqref{eq1.40}, an analogue of Theorem \ref{t1.5} is true in the sectors $\Omega_1$ and $\Omega_2$. Similarly, for $\varphi_1(\lambda,x)$ and $\varphi_2(\lambda,x)$, an analogue of Theorem \ref{t1.7} holds in the sectors $\Omega_1^-$ and $\Omega_2^-$ correspondingly.
\end{remark}

\begin{conclusion}
Theorems \ref{t1.5}, \ref{t1.7} provide simple formulas for the calculation of the function $m(x)$ on each half-axis $\mathbb{R}_\pm$ from the functions $\psi_0(\lambda,x)$ and $\varphi_0(\lambda,x)$. Therefore, one has to find the functions $\psi_0(\lambda,x)$ and $\varphi_0(\lambda,x)$ from the scattering data. Sections 2 and 3 are dedicated to these problems.
\end{conclusion}

\section{Transition matrix and direct scattering problem for the waves coming from $+\infty$}\label{s2}

{\bf 2.1} Expand every Jost solution $\{u_k(\lambda,x)\}_0^2$ by other Jost solutions $\{v_l(\lambda,x)\}_0^2$,
\begin{equation}
u_k(\lambda,x)=\sum\limits_{l=0}^2t_{k,l}(\lambda)v_l(\lambda,x)\quad(0\leq k\leq2),\label{eq2.1}
\end{equation}
or
\begin{equation}
u(\lambda,x)=T(\lambda)v(\lambda,x)\label{eq2.2}
\end{equation}
where $u(\lambda,x)\stackrel{\rm def}{=}\col[u_0(\lambda,x),u_1(\lambda,x),u_2(\lambda,x)]$, $v(\lambda,x)\stackrel{\rm def}{=}\col[v_0(\lambda,x),v_1(\lambda,x),v_2(\lambda,x)]$, and $T(\lambda)$ is the {\bf transition matrix},
\begin{equation}
T(\lambda)\stackrel{\rm def}{=}\left[
\begin{array}{ccc}
t_{00}(\lambda)&t_{0,1}(\lambda)&t_{0,2}(\lambda)\\
t_{1,0}(\lambda)&T_{1,1}(\lambda)&t_{1,2}(\lambda)\\
t_{2,0}(\lambda)&t_{2,1}(\lambda)&t_{2,2}(\lambda)
\end{array}\right].\label{eq2.3}
\end{equation}

\begin{remark}\label{r2.1}
It is sufficient to consider only one relation in \eqref{eq2.1}, e. g., for $k=0$,
\begin{equation}
u_0(\lambda,x)=t_{0,0}(\lambda)v_0(\lambda,x)+t_{0,1}(\lambda)v_1(\lambda,x)+t_{0,2}(\lambda)v_2(\lambda,x),\label{eq2.4}
\end{equation}
the rest of equalities in \eqref{eq2.1} follow from \eqref{eq2.4} after substitutions $\lambda\rightarrow\lambda\zeta_1$ and $\lambda\rightarrow\lambda\zeta_2$, due to \eqref{eq1.42} (Remark \ref{r1.4}). Besides,
\begin{equation}
\begin{array}{ccc}
t_{0,0}(\lambda\zeta_1)=t_{1,1}(\lambda);&t_{0,1}(\lambda\zeta_1)=t_{1,2}(\lambda):&t_{0,2}(\lambda\zeta_1)=t_{1,0}(\lambda);\\
t_{0,0}(\lambda\zeta_2)=t_{2,2}(\lambda);&t_{0,1}(\lambda\zeta_2)=t_{2,0}(\lambda);&t_{0,2}(\lambda\zeta_2)=t_{2,1}(\lambda).
\end{array}\label{eq2.5}
\end{equation}
Hence it follows that the transition matrix $T(\lambda)$ \eqref{eq2.3} is expressed via $t_{0,0}(\lambda)$, $t_{0,1}(\lambda)$, $t_{0,2}(\lambda)$ and
\begin{equation}
T(\lambda)=\left[
\begin{array}{ccc}
t_{0,0}(\lambda)&t_{0,1}(\lambda)&t_{0,2}(\lambda)\\
t_{0,2}(\lambda\zeta_1)&t_{0,0}(\lambda\zeta_1)&t_{0,1}(\lambda\zeta_1)\\
t_{0,1}(\lambda\zeta_2)&t_{0,2}(\lambda\zeta_2)&t_{0,0}(\lambda\zeta_2)
\end{array}\right]\label{eq2.6}
\end{equation}
\end{remark}

By $W_{k,s}(v,\lambda,x)$, we denote Vronskian of the functions $v_k(\lambda,x)$ and $v_s(\lambda,x)$,
\begin{equation}
\begin{array}{ccc}
W_{k,s}(v,\lambda,x)=W\{v_k(\lambda,x),v_s(\lambda,x)\}\stackrel{\rm def}{=}v_k(\lambda,x)v'_s(\lambda,x)-v_s(\lambda,x)v'_k(\lambda,x)\\
(0\leq k,s\leq2),
\end{array}\label{eq2.6}
\end{equation}
and let ``+'' be involution,
\begin{equation}
f^+(\lambda)=\overline{f(\overline{\lambda})}.\label{eq2.7}
\end{equation}

\begin{lemma}\label{l2.1}
For the Vronskians $W_{k,s}(v,\lambda,x)$ \eqref{eq2.6} the following representations hold:
\begin{equation}
\begin{array}{ccc}
W_{1,2}=\sqrt3z_+\zeta_0v_0^+(\lambda,x);\, W_{0,1}(v,\lambda,x)=\sqrt3z_+\zeta_2v_1^+(\lambda,x);\\
W_{2,0}(v,\lambda,x)=\sqrt3z_+\zeta_1v_2^+(\lambda,x)
\end{array}\label{eq2.8}
\end{equation}
($z_+=\lambda n_+$).
\end{lemma}

P r o o f. Since
$$W'_{k,s}(v,\lambda,x)=v_k(\lambda,x)v''_s(\lambda,x)-v_s(\lambda,x)v''_k(\lambda,x);$$
$$W''_{k,s}(v,\lambda,x)=v'_k(\lambda,x)v''_s(\lambda,x)-v'_s(\lambda,x)v''_k(\lambda,x),$$
then, upon differentiating the second equality one more time and using \eqref{eq1.14}, we obtain that $W_{k,s}(v,\lambda,x)$ is the solution to equation
\begin{equation}
iy'''(x)=-m(x)\lambda^3y(x)\label{eq2.9}
\end{equation}
which is obtained from \eqref{eq1.11} by applying the involution ``+'' \eqref{eq2.7}. To obtain \eqref{eq2.8}, it is necessary to take asymptotics (a) \eqref{eq1.16} into account. So, for $W_{1,2}(v,\lambda,x)$, we have
$$W_{1,2}(v,\lambda,x)=i\zeta_2z_+e^{iz_+(\zeta_1+\zeta_2)}-i\zeta_1z_+e^{iz_+(\zeta_1+\zeta_2)}=\sqrt3z_+e^{-iz_0z_+x}\,(x\rightarrow\infty).\blacksquare$$

For the Wronskians $W_{k,s}(u,\lambda,x)$,
\begin{equation}
W_{k,s}(u,\lambda,x)\stackrel{\rm def}{=}u_k(\lambda,x)u'_s(\lambda,x)-u_s(\lambda,x)u'_k(\lambda,x)\quad(0\leq k,s\leq2),\label{eq2.10}
\end{equation}
formulas analogous to \eqref{eq2.8} are true after the natural substitution $z_+\rightarrow z_-$ ($=\lambda n_-$).

Using \eqref{eq2.1}, calculate $W_{1,2}(u,\lambda,x)$,
$$W_{1,2}(u,\lambda,x)=[t_{1,0}(\lambda)v_0(\lambda,x)+t_{1,1}(\lambda)v_1(\lambda,x)+t_{1,2}(\lambda)v_2(\lambda,x)][t_{2,0}(\lambda)v_0(\lambda,x)$$
$$+t_{2,1}(\lambda)v_1(\lambda,x)+t_{2,2}(\lambda)v_2(\lambda,x)]'-[t_{2^0}v_0(\lambda,x)+t_{2^1}v_1(\lambda,x)+t_{2,2}(\lambda)v_2(\lambda,x)]$$
$$\times[t_{1,0}(\lambda)v_0(\lambda,x)+t_{1,1}(\lambda)v_1(\lambda,x)+t_{1,2}(\lambda)v_2(\lambda,x)]'=W_{0,1}(v,\lambda,x)\left[
\begin{array}{ccc}
t_{1,0}(\lambda)&t_{1,1}(\lambda)\\
t_{2,0}(\lambda)&t_{2,1}(\lambda)
\end{array}\right]$$
$$+W_{0,2}(v,\lambda,x)\left[
\begin{array}{ccc}
t_{1,0}(\lambda)&t_{1,2}(\lambda)\\
t_{2,0}(\lambda)&t_{2,2}(\lambda)
\end{array}\right]+W_{1,2}(v,\lambda,x)\left[
\begin{array}{ccc}
t_{1,1}(\lambda)&t_{1,2}(\lambda)\\
t_{2,1}(\lambda)&t_{2,2}(\lambda)
\end{array}\right]$$
and, taking into account \eqref{eq2.8}, we obtain
$$\sqrt3z_-u_0^+(\lambda,x)=\sqrt3z_+\{T_{0,0}(\lambda)v_0^+(\lambda,x)+\zeta_2T_{0,2}(\lambda)v_1^+(\lambda,x)+\zeta_1T_{0,1}(\lambda)v_2^+(\lambda,x)\}$$
where $T_{k,s}(\lambda)$ are the algebraic complements of the elements $t_{k,s}(\lambda)$ of the matrix $T(\lambda)$ \eqref{eq2.3}. As a result, we have
$$n_-u_0^+(\lambda,x)=n_+\{T_{0,0}(\lambda)v_0^+(\lambda,x)+\zeta_2T_{0,2}(\lambda)v_1^+(\lambda,x)+\zeta_1T_{0,1}(\lambda)v_2^+(\lambda,x)\}.$$
Similarly, calculating $W_{0,1}(u,\lambda,x)$ and $W_{0,2}(u,\lambda,x)$ and using \eqref{eq2.1} and \eqref{eq2.8}, we obtain
$$\zeta_2n_-u_1^+(\lambda,x)=n_+\{T_{1,0}(\lambda)v_0^+(\lambda,x)+\zeta_2T_{1,2}(\lambda)v_1^+(\lambda,x)+\zeta_1T_{0,1}(\lambda)v_2^+(\lambda,x)\};$$
$$\zeta_1n_-u_2^+(\lambda,x)=n_+\{T_{2,0}(\lambda)v_0^+(\lambda,x)+\zeta_2T_{2,2}(\lambda)v_1^+(\lambda,x)+\zeta_1T_{2,1}(\lambda)v_2^+(\lambda,x)\}.$$
Obtained equalities in the matrix form become
\begin{equation}
n_-Ju^+(\lambda,x)=n_+\widehat{T}(\lambda)Jv^+(\lambda,x),\label{eq2.11}
\end{equation}
where $u^+(\lambda,x)=(u(\lambda,x))^+$, $v^+(\lambda,x)=(v(\lambda,x))^+$, $\widehat{T}(\lambda)=[\widehat{T}_{k,s}(\lambda)]$ is the matrix of algebraic complements of $T(\lambda)$, and $J$ is an involution,
\begin{equation}
J=\left[
\begin{array}{ccc}
1&0&0\\
0&0&\zeta_2\\
0&\zeta_1&0
\end{array}\right]\quad(J^*=J,J^2=I).\label{eq2.12}
\end{equation}
Since ${\displaystyle T^{-1}(\lambda)=\frac1{\Delta(\lambda)}\{\widehat{T}^t(\lambda)\}}$, where $\Delta(\lambda)=\det T(\lambda)$ and $T^t(\lambda$ is the matrix obtained from $T(\lambda)$ upon transposition, then \eqref{eq2.11} implies
$$n_-J\{T^t(\lambda)\}Ju^+(\lambda,x)=n_+\Delta(\lambda)v^+(\lambda,x).$$
Applying ``+'' \eqref{eq2.7} to both sides of this equality, we obtain
$$n_-JT^*(\overline{\lambda})Ju(\lambda,x)=n_+\Delta^+(\lambda)v(\lambda,x).$$
Hence finding $v(\lambda,x)$ and substituting it into \eqref{eq2.2}, we have
$$\{n_-T(\lambda)JT^*(\overline{\lambda})J-n_+\Delta^+(\lambda)\}u(\lambda,x)$$
and thus
\begin{equation}
n_-T(\lambda)JT^*(\overline{\lambda})=n_+\Delta^+(\lambda)J.\label{eq2.13}
\end{equation}

\begin{lemma}\label{l2.2}
The transition matrix $T(\lambda)$ \eqref{eq2.3} has the properties
\begin{equation}
\begin{array}{lll}
({\rm i})\,{\displaystyle\det T(\lambda)=\frac{m_-}{m_+};}\\
({\rm ii})\,\displaystyle{T(\lambda)JT^*(\overline{\lambda})=k^2J\,\left(k\stackrel{\rm def}{=}\frac{n_-}{n_+}\right)}
\end{array}\label{eq2.14}
\end{equation}
\end{lemma}

P r o o f. Equation \eqref{eq2.13} implies that
\begin{equation}
m_-\Delta(\lambda)\Delta^+(\lambda)=m_+(\Delta^+(\lambda))^3.\label{eq2.15}
\end{equation}
Applying operation ``+'' \eqref{eq2.7} to this equation, we obtain $m_-\Delta^+(\lambda)\Delta(\lambda)=m_+(\Delta(\lambda))^3$. Subtracting this equality from \eqref{eq2.15}, we have $\Delta^3(\lambda)=(\Delta^+(\lambda))^3$, and thus $\Delta(\lambda)/\Delta^+(\lambda)=\zeta_p$ ($\zeta_p$ is one of the roots of \eqref{eq1.2}). Substituting $\Delta(\lambda)=\zeta_p\Delta^+(\lambda)$ into \eqref{eq2.15}, we obtain ${\displaystyle(\Delta^+(\lambda))^2\left(\Delta^+(\lambda)-\zeta_p\frac{m_-}{m_+}\right)=0.}$ Equality $\Delta^+(\lambda)=0$ contradicts linear independence of $\{u_k(\lambda,x)\}_0^2$. Therefore $\Delta^+(\lambda)=\zeta_pk^3$ ($k={\displaystyle\frac{n_-}{n_+}}$), substituting this expression for $\Delta^+(\lambda)$ into \eqref{eq2.13}, we have
$$T(\lambda)JT^*(\overline{\lambda})=k^2\zeta_pJ.$$
Equating the (1,1) elements in this equality, we obtain
\begin{equation}
k^2\zeta_p=t_{0,0}(\lambda)t_{0,0}^+(\lambda)+\zeta_1t_{0,1}(\lambda)t_{0,2}^+(\lambda)+\zeta_2t_{0,2}(\lambda)t_{0,1}^+(\lambda).\label{eq2.16}
\end{equation}
For $\lambda\in\mathbb{R}$, the right side of this equality is real and thus $\zeta_p=\zeta_0=1$, which concludes the proof. $\blacksquare$
\vspace{5mm}

{\bf 2.2} Turn to calculation of the coefficients $t_{0,0}(\lambda)$, $t_{0,1}(\lambda)$, $t_{0,1}(\lambda)$, $t_{0,2}(\lambda)$. Upon differentiating equality \eqref{eq2.4}, we obtain the system
\begin{equation}
\left\{
\begin{array}{lll}
u_0(\lambda,x)=t_{0,0}(\lambda)v_0(\lambda,x)+t_{0,1}(\lambda)v_1(\lambda,x)+t_{0,2}(\lambda)v_2(\lambda,x);\\
u'_0(\lambda,x)=t_{0,0}(\lambda)v'_0(\lambda,x)+t_{0,1}(\lambda)v'_1(\lambda,x)+t_{0,2}(\lambda)v'_2(\lambda,x);\\
u''_0(\lambda,x)=t_{0,0}(\lambda)v''_0(\lambda,x)+t_{0,1}(\lambda)v''_1(\lambda,x)+t_{0,2}(\lambda)v''_2(\lambda,x).
\end{array}\right.\label{eq2.17}
\end{equation}
In terms of the Wronskians \eqref{eq1.37}, the solution to system \eqref{eq2.17} is
\begin{equation}
\begin{array}{ccc}
t_{0,0}(\lambda)=W(u_0,v_1,v_2)/W(v_0,v_1,v_2);\quad t_{0,1}(\lambda)=W(v_0,u_0,v_2)/W(v_0,v_1,v_2);\\
t_{0,2}(\lambda)=W(v_0,v_1,u_0)/W(v_0,v_1,v_2).
\end{array}\label{eq2.18}
\end{equation}
Hence, due to \eqref{eq1.38} and \eqref{eq2.8}, follows the statement.

\begin{lemma}\label{l2.4}
The coefficients $\{t_{0,k}(\lambda)\}_0^2$ are given by
$$t_{0,0}(\lambda)=\frac{\zeta_0n_+}{3\lambda^2}\{u_0(\lambda,x)(v_0^+(\lambda,x))''-u'_0(\lambda,x)(v_0^+(\lambda,x))'+u''_0(\lambda,x)v_0^+(\lambda,x)\};$$
\begin{equation}
t_{0,1}(\lambda)=\frac{\zeta_1n_+}{3\lambda^2}\{u_0(\lambda,x)(v_2^+(\lambda,x))''-u'_0(\lambda,x)(v_2^+(\lambda,x))'+u''_0(\lambda,x)v_2^+(\lambda,x)\};\label{eq2.19}
\end{equation}
$$t_{0,2}(\lambda)=\frac{\zeta_2n_+}{3\lambda^2}\{u_0(\lambda,x)(v_1^+(\lambda,x))''-u'_0(\lambda,x)(v_1^+(\lambda,x))'+u''_0(\lambda,x)v_1^+(\lambda,x)\},$$
besides, the right sides of formulas \eqref{eq2.19} don't depend on $x$.
\end{lemma}

So, for $t_{0,0}(\lambda)$, e.g., we have
$$[u_0(\lambda,x)(v_0^+(\lambda,x))''-u'_0(\lambda,x)(v_0^+(\lambda,x)'+u''_0(\lambda,x)v_0^+(\lambda,x)]'$$
$$=u_0(\lambda,x)(v_0^+(\lambda,x))'''+u'''_0(\lambda,x)v_0^+(\lambda,x)=i\lambda^3m(x)u_0(\lambda,x)v_0^+(\lambda,x)$$
$$-i\lambda^3m(x)u_0(\lambda,x)v_0^+(\lambda,x)=0,$$
due to equations \eqref{eq1.11} and \eqref{eq2.9}.

Therefore, hereinafter we assume that in equalities \eqref{eq2.19} $x=0$.
\vspace{5mm}

{\bf 2.3} Study zeros of the function $t_{0,0}(\lambda)$ situated in $\Omega_0^-$. First of all, notice that $t_{0,0}(\lambda)$ is holomorphic in $\Omega_0^-$. If $w\in\Omega_0^-$ is a zero of the function $t_{0,0}(\lambda)$, $t_{0,0}(w)=0$, then \eqref{eq2.4} implies that
\begin{equation}
u_0(\lambda,x)=t_{0,1}(\lambda)v_1(\lambda,x)+t_{0,2}(\lambda)v_2(\lambda,x).\label{eq2.20}
\end{equation}
Reality of $w^3$ (see Remark \ref{r1.2}) yields that $w$ can lie on the rays $l_{\zeta_2}$, $\widehat{l}_{\zeta_1}$, i.e., either $w=\mu\zeta_2\in l_{\zeta_2}$ ($\mu>0$) or $w=\nu\zeta_1\in\widehat{l}_{\zeta_1}$ ($\nu<0$).

\begin{picture}(200,200)
\put(0,100){\vector(1,0){200}}
\put(100,0){\vector(0,1){200}}
\put(150,0){\vector(-1,2){100}}
\put(150,200){\vector(-1,-2){100}}
\put(200,67){\line(-3,1){200}}
\put(0,67){\line(3,1){200}}
\put(30,190){$l_{\zeta_1}$}
\put(135,150){$\nu\zeta_1$}
\put(80,170){$\Omega_1^-$}
\put(121,150){$\circ$}
\put(70,147){$\times$}
\put(150,96){$\times$}
\put(154,105){$\mu\zeta_0$}
\put(122,40){$\circ$}
\put(110,35){$\nu\zeta_1$}
\put(68,40){$\times$}
\put(80,35){$\mu\zeta_2$}
\put(60,5){$l_{\zeta_2}$}
\put(60,60){$\Omega_0^-$}
\put(50,96){$\circ$}
\put(38,105){$\nu\zeta_0$}
\put(190,105){$l_{\zeta_0}$}
\put(157,65){$\Omega_2^-$}
\put(50,150){$\mu\zeta_1$}
\qbezier(77,91)(100,70)(118,94)
\qbezier(100,120)(80,110)(80,93)
\qbezier(120,93)(120,120)(100,123)
\end{picture}

\hspace{20mm} Fig. 2

Points $\mu\zeta_0$ and $\nu\zeta_1$ are zeros of the function $t_{1,1}(\lambda)$ and lie in the sector $\Omega_2^-$. Analogously, $\mu\zeta_1$ and $\nu\zeta_0$, from the sector $\Omega_1^-$, are zeros of the function $t_{2,2}(\lambda)$. If $t_{0,1}(\mu\zeta_2)=0$, then \eqref{eq2.20} implies
\begin{equation}
u_0(\mu\zeta_2,x)=t_{0,2}(\mu\zeta_2)v_2(\mu\zeta_2,x).\label{eq2.21}
\end{equation}
Point $\mu\zeta_2$ belongs to two sectors $\Omega_0^-$ and $\Omega_2$ simultaneously. The left side of this equality belongs to $L^2(\mathbb{R}_-)$ (Theorem \ref{t1.6}) and the right, correspondingly, to $L^2(\mathbb{R}_+)$ (Theorem \ref{t1.4}). So, if $\mu\zeta_2$ is a zero of both functions $t_{0,0}(\lambda)$ and $t_{0,1}(\lambda)$, then the function $u_0(\mu\zeta_2,x)$ \eqref{eq2.21} belongs to $L^2(\mathbb{R})$ and generates a bound state of equation \eqref{eq1.11}.

Analogously, if $\nu\zeta_1\in\Omega_0^-\cap\Omega_1$ is the common zero for $t_{0,0}(\lambda)$ and $t_{0,2}(\lambda)$, then \eqref{eq2.4} yields
\begin{equation}
u_0(\nu\zeta_1,x)=t_{0,1}(\nu\zeta_1)v_1(\nu\zeta_1,x).\label{eq2.22}
\end{equation}

Theorems \ref{t1.4}, \ref{t1.6} imply that $u(\nu\zeta_1,x)\in L^2(\mathbb{R})$, and thus $u(\nu\zeta_1,x)$ \eqref{eq2.12} also is a bound state of equation \eqref{eq1.11}.

\begin{lemma}\label{l2.5}
Joint zeros $\mu\zeta_2\in l_{\zeta_2}$ of the functions $t_{0,0}(\lambda)$ and $t_{0,1}(\lambda)$ generate bound states $u(\mu\zeta_1,x)$ \eqref{eq2.21} of equation \eqref{eq1.11}. And joint zeros $\nu\zeta_1\in\widehat{l}_{\zeta_1}$ of the functions $t_{0,0}(\lambda)$ and $t_{0,2}(\lambda)$ also set bound states $u(\nu\zeta_1,x)$ \eqref{eq2.22} of equation \eqref{eq1.11}.
\end{lemma}

\begin{lemma}\label{l2.6}
The set of zeros $\mu\zeta_2$ (as well as $\nu\zeta_1$) from $\Omega_0^-$ of the function $t_{0,0}(\lambda)$ is no more than finite, moreover, all zeros are of multiplicity two.
\end{lemma}

P r o o f. Analiticity of the fuction $t_{0,0}(\lambda)$ in the sector $\Omega_0^-$ implies that in order to prove that the set of zeros $\{\mu\zeta_2\}$, $\{\nu\zeta_1\}$ is finite it is sufficient to prove that function $t_{0,0}(\lambda)$ does not vanish when $|\lambda|\gg1$ ($\lambda\in\Omega_0^-$). Equation \eqref{eq2.19} implies
\begin{equation}
t_{0,0}(\lambda)=\frac{n_+}{3\lambda^2}\{u_0(\lambda,0)(v_0^+(\lambda,0))''-u'_0(\lambda,0)(v_0^+(\lambda,0))'+u''_0(\lambda,0)v_0^+(\lambda,0)\}.\label{eq2.23}
\end{equation}
Taking into account \eqref{eq1.70} and \eqref{eq1.64}, we have
$$u_0(\lambda,0)(v_0^+(\lambda,0))''=-\lambda^2n_+^2\varphi_0(\lambda,0)\overline{\psi_0^2(\overline{\lambda},0)},$$
and analogously, from \eqref{eq1.71} and \eqref{eq1.65} we find that
$$u'_0(\lambda,0)(v_0^+(\lambda,0))'=\lambda^2n_+n_-\varphi_0^1(\lambda,0)\overline{\psi_0^1(\overline{\lambda},0)}.$$
Finally, using \eqref{eq1.72} and \eqref{eq1.62}, we obtain
$$u''_0(\lambda,0)v_0^+(\lambda,0)=-\lambda^2n_-^2\varphi_0^2(\lambda,0)\overline{\psi_0(\overline{\lambda},0)}.$$
As a result, $t_{0,0}(\lambda)$ \eqref{eq2.23} becomes
\begin{equation}
t_{0,0}(\lambda)=\frac{n_+}3\{-n_+^2\varphi_0(\lambda,0)\overline{\psi_0^2(\overline{\lambda},0)}+n_+n_-\varphi_0^1(\lambda,0)\overline{\psi_0^1(\overline{\lambda},0)}-n_-^2\varphi_0^2(\lambda,0)\overline{\psi_0(\overline{\lambda},0)}\}.\label{eq2.24}
\end{equation}
Let
$$r\stackrel{\rm def}{=}-n_+^2+n_+n_--n_-^2<0\quad(\forall n_+,n_->0)$$
and
$$p\stackrel{\rm def}{=}\frac{-n_+^2}r,\quad q\stackrel{\rm def}{=}\frac{-n_-^2}r,\quad l=\frac{n_+n_-}r<0,$$
besides,
$$p+q+l=1.$$
Using $l=1-p-q$, rewrite equality \eqref{eq2.24} as
$$t_{0,0}(\lambda)=\frac{n+r}3\{p\varphi_0(\lambda,0)\overline{\psi_0^2(\overline{\lambda},0)}+l\varphi_0^1(\lambda,0)\overline{\psi_0^1(\overline{\lambda},0)}+q\varphi_0^2(\lambda,0)\overline{\psi_0(\overline{\lambda},0)}\}$$
\begin{equation}
\begin{array}{ccc}
{\displaystyle=\frac{n+r}3\{p[\varphi_0(\lambda,0)\overline{\psi_0^2(\overline{\lambda},0)}-\varphi_0^1(\lambda,0)\overline{\psi_0^1(\overline{\lambda},0)}]+q[\varphi_0^2(\lambda,0)\overline{\psi_0(\overline{\lambda},0)}}\\
{\displaystyle-\varphi_0^1(\lambda,0)\overline{\psi_0^1(\overline{\lambda},0)}]+\varphi_0^1(\lambda,0)\overline{\psi_0^1(\overline{\lambda},0)}\}=\frac{n_+r}3\varphi_0^1(\lambda)\overline{\psi_0^1(\overline{\lambda},0)}}
\end{array}\label{eq2.25}
\end{equation}
$$\times\left\{p\left[\frac{\varphi_0(\lambda,0)\overline{\psi_0^2(\overline{\lambda},0)}}{\varphi_0^1(\lambda,0)\overline{\psi_0^1(\overline{\lambda},0)}}-1\right]+1+q\left[\frac{\varphi_0^2(\lambda,0)\overline{\psi_0(\overline{\lambda},0)}
}{\varphi_0^1(\lambda,0)
\overline{\psi_0^1(\overline{\lambda},0)}}-1\right]\right\}.$$
Taking into account the fact that for $\lambda\rightarrow\infty$ ($\lambda\in\Omega_0^-$), character of the growth of functions $\varphi_0(\lambda,0)$ and $\varphi_0^1(\lambda,0)$, $\varphi_0^2(\lambda,0)$ is the same and also that functions $\overline{\psi_0^1(\overline{\lambda},0)}$, $\overline{\psi_0^2(\overline{\lambda},0)}$, $\overline{\psi_0(\overline{\lambda},0)}$ behave in the same way when $\lambda\rightarrow\infty$ ($\lambda\in\Omega_0^-$), we obtain that expression in parenthesis in \eqref{eq2.25} tends to 1, $\{...\}\rightarrow1$ ($\lambda\rightarrow\infty$, $\lambda\in\Omega_0^-$). Notice that function $\varphi_0^1(\lambda,0)\overline{\psi_0^1(\overline{\lambda},0}$ for $\lambda=0$ equals 1, and for $\lambda\rightarrow\infty$ ($\lambda\in\Omega_0^-$) it increases in modulo. Hence we deduce that $t_{0,0}(\lambda)$ does not vanish when $|\lambda|\gg1$ ($\lambda\in\Omega_0^-$).

Now study multiplicity of zeros of the function $t_{0,0}(\lambda)$. Equation \eqref{eq2.19} implies
\begin{equation}
t_{0,0}(\lambda)=\frac{n_+}{3\lambda^2}g(\lambda,x)\label{eq2.26}
\end{equation}
where $g(\lambda,x)$ is given by
\begin{equation}
g(\lambda,x)\stackrel{\rm def}{=}u_0(\lambda,x)(v_0^+(\lambda,x))''-u'_0(\lambda,x)(v_0^+(\lambda,x))'+u''_0(\lambda,x)v_0^+(\lambda,x)\label{eq2.27}
\end{equation}
and does not depend on $x$ (Lemma \ref{l2.4}). Hereinafter, by the dot above we denote the derivative relative to $\lambda$, ${\displaystyle\dot{f}(\lambda)=\frac{df(\lambda)}{d\lambda}}$. Equationn \eqref{eq2.27} implies that
\begin{equation}
\begin{array}{ccc}
\dot{g}(\lambda,x)=\dot{u}_0(\lambda,x)(v_0^+(\lambda,x))''-\dot{u}'_0(\lambda,x)(v_0^+(\lambda,x))'+\dot{u}''_+(\lambda,x)v_0^+(\lambda,x)\\
+u_0(\lambda,x)(\dot{v}_0^+(\lambda,x))''-u_0(\lambda,x)(\dot{v}_0^+(\lambda,x))'+u''_0(\lambda,x)\dot{v}_0^+(\lambda,x).
\end{array}\label{eq2.28}
\end{equation}
Using equations
$$u'''_0(\lambda,x)=\frac{\lambda^3}im(x)u_0(\lambda,x);\quad\dot{u}'''_0(\lambda,x)=\frac{\lambda^3}im(x)\dot{u}_0(\lambda,x)+\frac{3\lambda^2}im(x)u_0(\lambda,x);$$
$$(v_0^+(\lambda,x))'''=-\frac{\lambda^3}im(x)v_0^+(\lambda,x);$$
$$(\dot{v}_0^+(\lambda,x))'''=-\frac{\lambda^3}im(x)\dot{v}_0^+(\lambda,x)-\frac{3\lambda^2}im(x)v_0^+(\lambda,x);$$
it is easy to show that $\dot{g}(\lambda,x)$ \eqref{eq2.28} does not depend on $x$, therefore, hereinafter we assume that $x=0$. Let $w$ be a zero of the function $g(\lambda,0)$, i.e.,
\begin{equation}
u(w,0)(v^+(w,0))''-u'_0(w,0)(v_0^+(w,0))'+u''_0(w,0)v_0^+(w,0)=0.\label{eq2.29}
\end{equation}
Equation \eqref{eq2.26} implies that
$$\dot{t}_{0,0}(\lambda)=\frac{n_+}3\frac{\dot{g}(\lambda,0)\lambda^2-2\lambda g(\lambda,0)}{\lambda^4},$$
therefore, due to \eqref{eq2.29}, $g(w,0)=0$,
$$t_{0,0}(w)=\frac{n_+}{3w^2}\dot{g}(w,0).$$
Prove that \eqref{eq2.29} yields $\dot{g}(w,0)=0$. Rewrite equality \eqref{eq2.28} as
$$\dot{g}(\lambda,0)=u_0(\lambda,0)(v_0^+(\lambda,0))''\cdot\dot{\ln}u_0(\lambda,0)({v}_0^+(\lambda,0))''-u'_0(\lambda,0)(v_0^+(\lambda,0))'$$
$$\times\dot{\ln}u'_0(\lambda,0)(v_0^+(\lambda,0))'+u''_0(\lambda,0)v_0^+(\lambda,0)\cdot\dot{\ln}u''_0(\lambda,0)v_0^+(\lambda,0).$$
For $\lambda=w$, \eqref{eq2.29} implies
$$u'_0(w,0)(v^+(w,0))'=u(w,0)(v_0^+(w,0))''+u''(w,0)v_0^+(w,0),$$
and thus
$$\dot{g}(w,0)=u(w,0)(v^+(w,0))''\dot{\ln}\left(\frac{u_0(w,0)(v_0^+(w,0))''}{u_0(w,0)(v_0^+(w,0))''+u''_0(w,0)v_0^+(w,0)}\right)$$
\begin{equation}
+u''(w,0)v^+(w,0)\cdot\dot{\ln}\left(\frac{u''_0(w,0)v_0^+(w,0)}{u_0(w,0)(v_0^+(w,0))''+u''_0(w,0)v_0^+(w,0)}\right)\label{eq2.30}
\end{equation}
$$=u''_0(w,0)v_0^+(w,0)\left\{f(w)\dot{\ln}\left(\frac{f(w)}{1+f(w)}\right)+\dot{\ln}\left(\frac1{1+f(w)}\right)\right\}$$
where
$$f(w)=\frac{u(w,0)(v_0^+(w,0))''}{u''_0(w,0)v_0^+(w,0)}.$$
Simplify expression in the parenthesis; let
$$\ln\frac{f(w)}{1+f(w)}=F(w);\quad\frac{f(w)}{1+f(w)}=e^{F(w)},$$
and thus
$$f(w)=\frac{e^{F(w)}}{1-e^{F(w)}};\quad1+f(w)=\frac1{1-e^{F(w)}}.$$
Hence it follows that expression in parenthesis in \eqref{eq2.30} equals
$$\{...\}=\frac{e^{F(w)}}{1-e^{F(w)}}\cdot\dot{F}(w)+\dot{\ln}(1-F(w))=\frac{e^{F(w)}}{1-e^{F(w)}}\dot{F}(w)-\frac{e^{F(w)}}{1-e^{F(w)}}\cdot\dot{F}(w)$$
$$=0,$$
q.e.d. $\blacksquare$

\begin{remark}\label{r2.2}
Equation \eqref{eq2.5} implies that $t_{1,1}(\lambda)$ and $t_{2,2}(\lambda)$ in sectors $\Omega_2^-$ and $\Omega_1^-$ correspondingly have finite number of zeros and all of them are of multiplicity $2$.
\end{remark}
\vspace{5mm}

{\bf 2.4} Rewrite equality \eqref{eq2.4} as
\begin{equation}
r_0(\lambda)u_0(\lambda,x)=v_0(\lambda,x)+s_1(\lambda)v_1(\lambda,x)+s_2(\lambda)v_2(\lambda,x)\label{eq2.31}
\end{equation}
where
\begin{equation}
r_0(\lambda)\stackrel{\rm def}{=}\frac1{t_{0,0}(\lambda)};\quad s_1(\lambda)\stackrel{\rm def}{=}\frac{t_{0,1}(\lambda)}{t_{0,0}(\lambda)};\quad s_2(\lambda)=\frac{t_{0,2}(\lambda)}{t_{0,0}(\lambda)}.\label{eq2.32}
\end{equation}

\begin{remark}\label{r2.3}
For $k\rightarrow+\infty$, function $r_0(\lambda)u_0(\lambda,x)$ has asymptotic
$$r_0(\lambda)u_0(\lambda,x)\rightarrow e^{iz_+\zeta_0x}+s_1(\lambda)e^{iz_+\zeta_1x}+s_2(\lambda)e^{iz_+\zeta_2x}\quad(z_+=\lambda n_+),$$
and thus the incident (from $+\infty$) wave $e^{iz_+\zeta_0x}$ has the {\bf reflected (scattered)} wave $s_1(\lambda)e^{iz_+\zeta_1x}+s_2(\lambda)e^{iz_+\zeta_2x}$ where $s_1(\lambda)$ and $s_2(\lambda)$ are scattering coefficients of the incident wave $e^{iz_+\zeta_0x}$.

Asympotic behavior of $r_0(\lambda)u_0(\lambda,x)$ \eqref{eq2.31}, for $x\rightarrow-\infty$, is
$$r_0(\lambda)u_0(\lambda,x)\rightarrow r_0(\lambda)e^{iz_-\zeta_0x}\quad(z_-=\lambda n_-),$$
therefore, it is natural to consider $r(\lambda)$ as a {\bf coefficient of transformation} of wave $e^{iz_+\zeta_0x}$ into wave $e^{iz\zeta_0x}$.

Equality \eqref{eq2.16}, written in terms of $r_0(\lambda)$, $s_1(\lambda)$, $s_2(\lambda)$ \eqref{eq2.32}, becomes
\begin{equation}
k^2r_0(\lambda)r_0^+(\lambda)=1+\zeta_1s_1(\lambda)s_2^+(\lambda)+\zeta_2s_2(\lambda)s_1^+(\lambda),\label{eq2.33}
\end{equation}
this is unitarity property for the scattering problem for incident wave $e^{iz_+\zeta_0x}$.
\end{remark}

Upon substituting $\lambda\rightarrow\lambda\zeta_1$ and $\lambda\rightarrow\lambda\zeta_2$ into equality \eqref{eq2.31}, we obtain
\begin{equation}\begin{array}{ccc}
({\rm i})\quad r_1(\lambda)u_1(\lambda,x)=s_2(\lambda\zeta_1)u_0(\lambda,x)+v_1(\lambda,x)+s_1(\lambda\zeta_1)v_2(\lambda,x);\\
({\rm ii})\quad r_2(\lambda)u_2(\lambda,x)=s_1(\lambda\zeta_2)v_0(\lambda,x)+s_2(\lambda\zeta_2)v_1(\lambda,x)+v_2(\lambda,x).
\end{array}\label{eq2.33}
\end{equation}
Equalities \eqref{eq2.33} describe scattering of the incident (from $+\infty$) waves $e^{iz_+\zeta_1x}$ and $e^{iz_+\zeta_2x}$ ($r_k(\lambda)=r_0(\lambda\zeta_k)$, $k=1$, 2).

Since
$$T(\lambda)=R^{-1}(\lambda)(I+S(\lambda)),$$
where
\begin{equation}
R(\lambda)\stackrel{\rm def}{=}\left[
\begin{array}{ccc}
r_0(\lambda)&0&0\\
0&r_1(\lambda)&0\\
0&0&r_2(\lambda)
\end{array}\right];\quad S(\lambda)\stackrel{\rm def}{=}\left[
\begin{array}{ccc}
0&s_1(\lambda)&s_2(\lambda)\\
s_2(\lambda\zeta_1)&0&s_1(\lambda\zeta_1)\\
s_1(\lambda\zeta_2)&s_2(\lambda\zeta_2)&0
\end{array}\right],\label{eq2.34}
\end{equation}
equation \eqref{eq2.14} implies
\begin{equation}
[I+S(\lambda)]J[I+S^+(\overline{\lambda})]=k^2R(\lambda)JR^+(\overline{\lambda}).\label{eq2.35}
\end{equation}
Matrix $R(\lambda)$ \eqref{eq2.34} is said to be the {\bf transformation matrix}, and $S(\lambda)$ \eqref{eq2.34}, correspondingly, the {\bf scattering matrix}. Equality \eqref{eq2.35} is the {\bf law of conservation of energetic balance} between scattered and transformed waves.
\vspace{5mm}

{\bf 2.5} Calculate the Wronskians $W\{u_0(\lambda,x),v_k(\lambda,x)\}$ ($k=1$ 2), where $u_0(\lambda)$ is given by \eqref{eq2.31}, then, taking into account \eqref{eq2.8}, we obtain
\begin{equation}
\begin{array}{lll}
({\rm i})& r_0(\lambda)f_{0,1}(\lambda,x)=\sqrt3z_+\zeta_2v_1^+(\lambda,x)-s_2(\lambda)\sqrt3z_+\zeta_0v_0^+(\lambda,x)&(il_{\zeta_2});\\
({\rm ii})& r_0(\lambda)f_{0,2}(\lambda,x)=-\sqrt3z_+\zeta_1v_2^+(\lambda,x)+s_1(\lambda)\sqrt3z_+\zeta_0v_0^+(\lambda,x)&(il_{\zeta_1}),
\end{array}\label{eq2.36}
\end{equation}
where
\begin{equation}
f_{0,k}(\lambda,x)=W\{u_0(\lambda,x),v_k(\lambda,x)\}\quad(k=1,2).\label{eq2.37}
\end{equation}

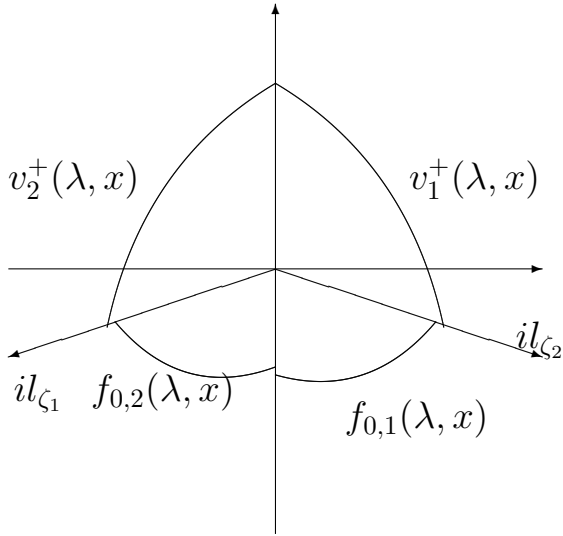
\begin{picture}(400,200)
\put(24,100){\vector(1,0){200}}
\put(124,0){\vector(0,1){200}}
\put(124,100){\vector(3,-1){100}}
\put(124,100){\vector(-3,-1){100}}
\put(26,50){$il_{\zeta_1}$}
\put(150,40){$f_{0,1}(\lambda,x)$}
\put(214,70){$il_{\zeta_2}$}
\qbezier(124,60)(159,50)(184,80)
\qbezier(124,63)(89,50)(64,80)
\qbezier(124,170)(74,140)(61,78)
\qbezier(124,170)(174,140)(187,78)
\put(174,130){$v_1^+(\lambda,x)$}
\put(24,130){$v_2^+(\lambda,x)$}
\put(55,50){$f_{0,2}(\lambda,x)$}
\end{picture}

\hspace{20mm} Fig. 3

Functions $v_1^+(\lambda,x)$ and $v_2^+(\lambda,x)$ are holomorphic in the sectors $\Omega_2^-$ and $\Omega_1^-$, and $f_{0,k}(\lambda,x)$ are analytic in the sectors $\Omega_0^-\cap\Omega_k$ ($k=1$, 2) correspondingly. Equalities \eqref{eq2.36} generate two jump problems on the rays $il_{\zeta_1}$ and $il_{\zeta_2}$ (Fig. 3).

Analogously, calculating Wronskians $W\{u_1(\lambda,x),v_k(\lambda,x)\}$ ($k=0$, 2) and using (i) \eqref{eq2.33} and \eqref{eq2.8}, we have
\begin{equation}
\begin{array}{lll}
({\rm i})&r_1(\lambda)f_{1,0}(\lambda,x)=-\sqrt3z_+\zeta_2v_1^+(\lambda,x)+s_1(\lambda\zeta_1)\sqrt3z_+\zeta_1v_2^+(\lambda,x)&(il_{\zeta_0});\\
({\rm ii})&r_1(\lambda)f_{1,2}(\lambda,x)=-s_2(\lambda\zeta_1)z_+\zeta_1v_2^+(\lambda,x)+\sqrt3z_+\zeta_0v_0^+(\lambda,x)&(il_{\zeta_1}),
\end{array}\label{eq2.38}
\end{equation}
where
\begin{equation}
f_{1,k}(\lambda,x)=W\{u_1(\lambda,x),v_k(\lambda,x)\}\quad(k=0,2).\label{eq2.39}
\end{equation}

\begin{picture}(300,200)
\put(70,100){\vector(1,0){200}}
\put(170,0){\vector(0,1){200}}
\put(70,133){\vector(3,-1){200}}
\put(170,100){\vector(-3,-1){100}}
\put(185,150){$v_1^+(\lambda,x)$}
\put(172,190){$il_{\zeta_0}$}
\put(190,50){$v_0^+(\lambda,x)$}
\put(76,85){$f_{1,2}(\lambda,x)$}
\put(72,57){$il_{\zeta_1}$}
\put(104,125){$f_{1,0}(\lambda,x)$}
\qbezier(193,89)(206,140)(170,133)
\qbezier(203,90)(170,50)(139,88)
\qbezier(149,108)(130,105)(139,88)
\qbezier(150,107)(150,120)(170,130)
\end{picture}

\hspace{20mm} Fig. 4.

Functions $v_0^+(\lambda,x)$ and $v_2^+(\lambda,x)$ are analytic in the sectors $\Omega_0^-$ and $\Omega_2^-$ and the functions $f_{1,k}(\lambda,x)$ in $\Omega_-^-\cap\Omega_k$ ($k=0$, 2) correspondingly. Thus, we have two jump problems on the rays $il_{\zeta_0}$ and $il_{\zeta_1}$. Finally, using (ii) \eqref{eq2.33}, calculate the Wronskians $W\{u_2(\lambda,x),v_k(\lambda,x)\}$ ($k=0$, 1), then, taking into account \eqref{eq2.8}, we obtain
\begin{equation}
\begin{array}{lll}
({\rm i})&r_2(\lambda)f_{2,0}(\lambda,x)=-s_2(\lambda\zeta_2)\sqrt3\zeta_2v_1^+(\lambda,x)+\sqrt3z_+\zeta_1v_2^+(\lambda,x)&(il_{\zeta_2});\\
({\rm ii})&r_2(\lambda)f_{2,1}(\lambda,x)=-s_1(\lambda\zeta_2)\sqrt3\zeta_2v_2^+(\lambda,x)+\sqrt3z_+v_0^+(\lambda,x)&(il_{\zeta_1})
\end{array}\label{eq2.40}
\end{equation}
where
\begin{equation}
f_{2,k}(\lambda,x)=W\{u_2(\lambda,x),v_k(\lambda,x)\}\quad(k=0,1).\label{eq2.41}
\end{equation}

\begin{picture}(300,200)
\put(70,100){\vector(1,0){200}}
\put(170,0){\vector(0,1){200}}
\put(260,131){\vector(-3,-1){200}}
\put(170,100){\vector(3,-1){100}}
\put(172,190){$il_{\zeta_0}$}
\put(200,140){$f_{2,0}(\lambda,x)$}
\put(190,50){$v_0^+(\lambda,x)$}
\put(260,120){$f_{2,1}(\lambda,x)$}
\put(104,125){$v_2^+(\lambda,x)$}
\put(247,58){$il_{\zeta_2}$}
\qbezier(193,110)(206,140)(170,133)
\qbezier(195,110)(220,120)(202,91)
\qbezier(203,88)(170,50)(139,88)
\qbezier(135,93)(130,120)(170,130)
\end{picture}

\hspace{20mm} Fig. 5.

Functions $v_0^+(\lambda,x)$ and $v_2^+(\lambda,x)$ are holomorphic in the sectors $\Omega_0^-$ and $\Omega_1^-$, and $f_{2,k}(\lambda,x)$, correspondingly, in the sectors $\Omega_1^-\cap\Omega_k$ ($k=0$, 1). Thus we obtain two more jump problems on the rays $il_{\zeta_0}$ and $il_{\zeta_2}$ (see Fig. 5).

\begin{remark}\label{r2.4}
Equalities \eqref{eq2.38} and \eqref{eq2.40} follow from \eqref{eq2.36} upon substituting $\lambda\rightarrow\lambda\zeta_1$ and $\lambda\rightarrow\lambda\zeta_2$.
\end{remark}

Define the functions holomorphic in respective sectors:
$$\psi_{2,1}(\lambda,x)\stackrel{\rm def}{=}r_2(\lambda)f_{2,1}(\lambda,x)e^{iz_+\zeta_0x}\quad(\lambda\in\Omega_2^-\cap\Omega_1);$$
\begin{equation}
\begin{array}{lll}
\psi_{1,2}(\lambda,x)\stackrel{\rm def}{=}r_1(\lambda)f_{1,2}(\lambda,x)e^{iz_+\zeta_0x}&(\lambda\in\Omega_1^-\cap\Omega_2);\\
\psi_{2,0}(\lambda,x)\stackrel{\rm def}{=}r_2(\lambda)f_{2,0}(\lambda,x)e^{iz_+\zeta_1x}&(\lambda\in\Omega_2^-\cap\Omega_0);
\end{array}\label{eq2.42}
\end{equation}
$$\psi_{1,0}(\lambda,x)\stackrel{\rm def}{=}r_1(\lambda)f_{1,0}(\lambda,x)e^{iz_+\zeta_2x}\quad(\lambda\in\Omega_1^-\cap\Omega_0).$$
Equations (ii) \eqref{eq2.31} and (ii) \eqref{eq2.40} yield two jump problems on the rays $il_{\zeta_1}$ and $il_{\zeta_2}$,
\begin{equation}
\begin{array}{lll}
\psi_{1,2}(\lambda,x)-\sqrt3\zeta_0\psi_0^{-1}(\lambda,x)=p_2(\lambda,x)v_2^+(\lambda,x)&(il_{\zeta_1});\\
\psi_{2,1}(\lambda,x)-\sqrt3\zeta_0\psi_0^+(\lambda,x)=p_1(\lambda,x)v_1^+(\lambda,x)&(il_{\zeta_2})\end{array}\label{eq2.43}
\end{equation}
where
\begin{equation}
p_2(\lambda,x)\stackrel{\rm def}{=}s_1(\lambda\zeta_1)\sqrt3\zeta_1e^{iz_+\zeta_0x};\quad p_1(\lambda,x)\stackrel{\rm def}{=}s_2(\lambda\zeta_2)\sqrt3\zeta_2e^{iz_+\zeta_0x}.\label{eq2.44}
\end{equation}
Boundary values of function $\psi_{2,2}(\lambda,x)$ on the ray $i\widehat{l}_{\zeta_1}$ coincide with boundary values of function $\psi_{1,0}(\lambda,x)$ on the ray $i\widehat{l}_{\zeta_2}$ since
$$\left.\psi_{2,1}(\lambda,x)\right|_{\lambda=-i\tau\zeta_1\in i\widehat{l}_{\zeta_1}}=r_0(-i\tau)f_{0,2}(-i\tau,x)e^{\tau n_+\zeta_1x};$$
$$\left.\psi_{1,0}(\lambda,x)\right|_{\lambda=-i\tau\zeta_2\in\widehat{l}_{\zeta_2}}=r_0(-i\tau)f_{0,2}(-i\tau,\lambda)e^{\tau n_+\zeta_1x}.$$

\begin{picture}(200,200)
\put(0,100){\vector(1,0){200}}
\put(100,0){\vector(0,1){200}}
\put(100,100){\line(-3,1){100}}
\put(100,100){\line(3,1){100}}
\put(2,134){$i\widehat{l}_{\zeta_2}$}
\put(25,144){$\psi_{1,0}(\lambda,x)$}
\put(100,150){$\psi_{2,0}(\lambda,x)$}
\qbezier(100,140)(135,145)(141,115)
\qbezier(100,135)(135,145)(138,112)
\qbezier(140,100)(155,105)(143,116)
\qbezier(56,100)(50,100)(58,115)
\qbezier(60,100)(58,105)(63,112)
\qbezier(100,141)(65,150)(57,116)
\put(0,105){$\psi_{1,2}(\lambda,x)$}
\put(170,110){$\psi_{2,1}(\lambda,x)$}
\put(175,135){$il_{\zeta_1}$}
\end{picture}

\hspace{20mm} Fig. 6.

Similarly, it is proved that boundary values of $\psi_{1,2}(\lambda,x)$ on $il_{\zeta_2}$ coincide with boundary values of $\psi_{2,0}(\lambda,x)$ on the ray $i\widehat{l}_{\zeta_1}$.

Define the function $\widehat{\psi}_{1,0}(\lambda,x)\stackrel{\rm def}{=}\psi_{1,0}(\widehat{\lambda},x)$ where $\widehat{\lambda}$ is a point from the sector $\Omega_0\cap\Omega_1^-$ symmetrical to $\lambda$ relative to the ray $il_{\zeta_0}$. Analogously, we define function $\widehat{\psi}_{2,0}(\lambda,x)$ relative to symmetry in the sector $\Omega_0\cap\Omega_1^-$. Since functions $\psi_{2,1}(\lambda,x)$ and $\widehat{\psi}_{1,0}(\lambda,x)$ coincide on the ray $i\widehat{l}_{\zeta_1}$, $\widehat{\psi}_{1,0}(\lambda,x)$ is an analytic extension of $\psi_{2,1}(\lambda,x)$ into sector $\Omega_2$. Analogously, $\widehat{\psi}_{2,0}(\lambda,x)$  is a holomorphic extension of function $\psi_{1,2}(\lambda,x)$ into sector $\Omega_1^-$. To obtain the jump problem on the ray $il_{\zeta_0},$, we multiply equality (i) \eqref{eq2.38} by $e^{iz_+\zeta_1x}$ and equality (i) \eqref{eq2.40} by $e^{iz_+\zeta_2x}$ and subtract the obtained, then we have
\begin{equation}
\psi_{2,0}(\lambda,x)-\psi_{1,0}(\lambda,x)=v_1^+(\lambda,x)p_3(\lambda,x)+v_2^+(\lambda,x)p_4(\lambda,x)\quad(il_{\zeta_0})\label{eq2.45}
\end{equation}
where
\begin{equation}
p_3(\lambda,x)\stackrel{\rm def}{=}\sqrt3z_+\zeta_2[s_2(\lambda\zeta_2)e^{iz_+\zeta_1x}-e^{iz_+\zeta_2x}];\, p_4(\lambda,x)\stackrel{\rm def}{=}\sqrt3z_+\zeta_1[e^{iz_+\zeta_1x}-s_1(\lambda\zeta_1)e^{iz_+\zeta_2x}].\label{eq2.46}
\end{equation}

\begin{picture}(200,200)
\put(0,100){\line(1,0){200}}
\put(100,0){\vector(0,1){200}}
\put(0,133){\vector(3,-1){200}}
\put(200,133){\vector(-3,-1){200}}
\put(100,150){$\widehat{\psi}_{1,0}(\lambda,x)$}
\put(100,185){$il_{\zeta_0}$}
\qbezier(100,140)(135,145)(141,115)
\qbezier(140,85)(155,105)(143,116)
\qbezier(100,60)(135,70)(144,85)
\qbezier(100,58)(70,60)(57,82)
\qbezier(56,86)(50,100)(58,115)
\qbezier(100,141)(65,150)(57,116)
\put(10,55){$il_{\zeta_1}$}
\put(0,105){$\psi_{1,2}(\lambda,x)$}
\put(170,110){$\psi_{2,1}(\lambda,x)$}
\put(173,53){$il_{\zeta_2}$}
\put(105,40){$\psi_0^+(\lambda,x)$}
\end{picture}

\hspace{20mm} Fig. 7.

So, in every sector $\Omega_0^-$, $\Omega_1^-$, $\Omega_2^-$ we have holomorphic functions (see Fig. 7) that on the rays $il_{\zeta_1}$, $il_{\zeta_2}$, $il_{\zeta_0}$ satisfy jump problem \eqref{eq2.43} and \eqref{eq2.45}. Inside sector $\Omega_0^-$, we have a holomorphic function $\psi_0^+(\lambda,x)$, and in the sectors $\Omega_1\cap\Omega_2^-$, $\Omega_0\cap\Omega_2^-$, $\Omega_0\cap\Omega_1^-$, $\Omega_2\cap\Omega_1^-$, correspondingly, $\psi_{2,1}(\lambda,x)$, $\widehat{\psi}_{1,0}(\lambda,x)$, $\widehat{\psi}_{2,0}(\lambda,x)$, $\psi_{1,2}(\lambda,x)$ with poles at the points $\{\mu_n\}\in l_{\zeta_0}$, $\{-\mu_n\zeta_2\}\in\widehat{l}_{\zeta_2}$, $\{-\nu_m\zeta_1\}\in l_{\zeta_1}$, $\{\nu_m\}\in\widehat{l}_{\zeta_0}$.

Consider a piecewise holomorphic in the sectors function (see Fig. 7),
\begin{equation}
F(\lambda,x)\stackrel{\rm def}{=}\left\{
\begin{array}{lllll}
\psi_0^+(\lambda,x)&(\lambda\in\Omega_0^-);\\
\psi_{2,1}(\lambda,x)&(\lambda\in\Omega_1\cap\Omega_2^-);\\
\widehat{\psi}_{1,0}(\lambda,x)&(\lambda\in\Omega_0\cap\Omega_2^-);\\
\psi_{2,0}(\lambda,x)&(\lambda\in\Omega_0\cap\Omega_1^-);\\
\psi_{1,2}(\lambda,x)&(\lambda\in\Omega_2\cap\Omega_1^-).
\end{array}\right.\label{eq2.47}
\end{equation}
And let
\begin{equation}
F_Q(\lambda,x)=F(\lambda,x)/Q(\lambda)\label{eq2.48}
\end{equation}
where $Q(\lambda)$ is an entire function not having zeros and growing, when $\lambda\rightarrow\infty$, faster than $\exp|z_+|(\sigma_+(0)+\sigma_-(0))$ (e.g., $Q(\lambda)=\exp z_+^2(\sigma_+(0)+\sigma_-(0))$). Then  \eqref{eq1.61} and \eqref{eq1.69} imply that $F_Q(\lambda,x)\rightarrow0$, for $\lambda\rightarrow\infty$, in every sector $\{\Omega_k^-\}_0^2$. In this case, $F_Q(\lambda,x)$ \eqref{eq2.48} is determined from its jumps on the rays $\{il_{\zeta_k}\}_0^2$ via a Cauchy type integral taking into account poles in corresponding sectors \cite{16,17},
$$F_Q(\lambda,x)=\sum\limits_n\frac{R_n(\zeta_0,x)}{(\lambda-\mu_n)^2}+\sum\limits_n\frac{R_n(\zeta_2,x)}{(\lambda+\mu_n\zeta_2)^2}+\sum\limits_m\frac{\widehat{R}_m(\zeta_0,x)}{(\lambda-\nu_m)^2}+\sum\limits_m\frac{\widehat{R}_m(\zeta_1,x)}{(\lambda+
\nu_m\zeta_1)^2}$$
\begin{equation}
+\frac1{2\pi i}\int\limits_0^\infty\frac{p_2(i\tau\zeta_1,x)v_1^+(i\tau,x)}{Q(i\tau\zeta_1)}\frac{d\tau}{\tau+i\zeta_2\lambda}+\frac1{2\pi i}\int\limits_0^\infty\frac{p_1(i\tau\zeta_2,x)v_2^+(i\tau,x)}{Q(i\tau\zeta_2)}\frac{d\tau}{\tau+i\zeta_1\lambda}\label{eq2.49}
\end{equation}
$$+\frac1{2\pi i}\int\limits_0^\infty\frac{p_3(i\tau,x)v_1^+(i\tau,x)+p_1(i\tau,x)v_2(i\tau,x)}{Q(i\tau)}\frac{d\tau}{\tau+i\lambda}.$$
Using relations $\left.g_{2,1}(\lambda,x)\right|_{\lambda=\lambda\zeta_2}=g_{1,0}(\lambda,x)$ and $\left.g_{1,2}(\lambda,x)\right|_{\lambda=\lambda\zeta_1}=g_{2,0}(\lambda,x)$, it is easy to show that
$$R_n(\zeta_2,x)=R_n(\zeta_0)\zeta_1e^{i(\zeta_0-\zeta_2)n_+\mu_nx};\quad\widehat{R}_m(\zeta_1,x)=\widehat{R}_m(\zeta_0,x)\zeta_2e^{i(\zeta_0-\zeta_1)n_+\nu_mx}.$$

Therefore, equality \eqref{eq2.49} becomes
$$F_Q(\lambda,x)=\sum\limits_nR_n(\lambda,x)\left\{\frac1{(\lambda-\mu_n)^2}+\frac{\zeta_1e^{i(\zeta_0-\zeta_2)n_+\mu_nx}}{(\lambda+\mu_n\zeta_2)^2}\right\}+\sum\limits_m\widehat{R}_m(\zeta_0,x)\left\{\frac1{(\lambda-\nu_m)^2}\right.$$
\begin{equation}
\left.+\frac{\zeta_2e^{i(\zeta_0-\zeta_1)n_+\nu_mx}}{(\lambda+\nu_m\zeta_1)^2}\right\}+\frac1{2\pi i}\int\limits_0^\infty\frac{p_2(i\tau\zeta_1,x)v_1^+(i\tau,x)}{Q(i\tau\zeta_1)}\frac{d\tau}{\tau+i\zeta_2\lambda}\label{eq2.50}
\end{equation}
$$+\frac1{2\pi i}\int\limits_0^\infty\frac{p_1(i\tau\zeta_2,x)v_2^+(i\tau,x)}{Q(i\tau\zeta_2)}\frac{d\tau}{\tau+i\zeta_1\lambda}$$
$$+\frac1{2\pi i}\int\limits_0^\infty\frac{p_3(i\tau,x)v_1^+(i\tau,x)+p_4(i\tau,x)v_2^+(i\tau,x)}{Q(i\tau)}\frac{d\tau}{\tau+i\lambda}.$$
For $\lambda\in\Omega_0^-$, $F_Q(\lambda,x)=v_0^+(\lambda,x)e^{iz_+\zeta_0x}/B(\lambda)$, hence, calculating the boundary values on the rays $\lambda\rightarrow it\zeta_1\in l_{\zeta_1}$ and $\lambda\rightarrow it\zeta_2\in l_{\zeta_2}$, and also using Sokhotski formulas \cite{16,17}, we obtain
$$({\rm i})\,v_2^+(it,x)e^{-n_+t\zeta_1x}Q^{-1}(it\zeta_1)=\sum\limits_nR_n(\zeta_0,x)\left[\frac1{(it\zeta_1-\mu_n)^2}+\frac{\zeta_1e^{i(\zeta_0-\zeta_2)n_+\mu_nx}}{(it\zeta_1+\mu_n\zeta_2)^2}\right]$$
$$+\sum\limits_m\widehat{R}_m(\zeta_0,x)\left[\frac1{(it\zeta_1-\nu_m)^2}+\frac{\zeta_2e^{i(\zeta_0-\zeta_1)n_+\nu_mx}}{(it\zeta_1-\nu_m\zeta_1)^2}\right]+\frac12\frac{p_2(it\zeta_1,x)v_1^+(it,x)}{Q(it\zeta_1)}$$
$$+\frac1{2\pi i}\int\limits_0^\infty\hspace{-4mm}/\frac{p_2(i\tau\zeta_1,x)v_1^+(i\tau,x)}{Q(i\tau\zeta_1)}\frac{d\tau}{t-\tau}+\frac1{2\pi i}\int\limits_0^\infty\frac{p_1(i\tau\zeta_2,x)v_2^+(i\tau,x)}{Q(i\tau\zeta_2)}\frac{d\tau}{\tau-\zeta_2t}$$
\begin{equation}
+\frac1{2\pi i}\int\limits_0^\infty\frac{p_3(it,x)v_1^+(it,x)+p_4(it,x)v_2^+(it,x)}{Q(it)}\frac{d\tau}{\tau-\zeta_1t};\label{eq2.51}
\end{equation}
$$({\rm ii})\,v_1^+(it,x)e^{-n_+\zeta_2tx}B(it\zeta_2)=\sum\limits_nR_n(\zeta_0x)\left[\frac1{(it\zeta_2-\mu_n)^2}+\frac{\zeta_1e^{i(\zeta_0-\zeta_1)n_+\mu_nx}}{(it\zeta_2+\mu_n\zeta_2)^2}\right]$$
$$+\sum\limits_mR_m(\zeta_0,x)\left[\frac1{(it\zeta_2-\nu_m)^2}+\frac{\zeta_2e^{i(\zeta_0-\zeta_2)n_+\nu_mx}}{(it\zeta_2-\nu_m\zeta_1)^2}\right]+\frac1{2\pi i}\int\limits_0^\infty\frac{p_2(it\zeta_1,x)v_1^+(it,x)}{Q(i\tau\zeta_1)}\frac{d\tau}{\tau-\zeta_1t}$$
$$+\frac12\frac{p_2(it\zeta_1)v_1^+(it,x)}{Q(it\zeta_1)}+\frac1{2\pi i}\int\limits_0^\infty\hspace{-4mm}/\frac{p_1(i\tau\zeta_2x)v_2^+(i\tau,x)}{Q(i\tau\zeta_2)}\frac{d\tau}{\tau-t}$$
$$+\frac1{2\pi i}\int\limits_0^\infty\frac{p_3(i\tau,x)v_1^+(i\tau,x)+p_4(i\tau,x)v_2^+(i\tau,x)}{Q(i\tau)}\frac{d\tau}{\tau-\zeta_2t}.$$

Multiply equality \eqref{eq2.50} by $(\lambda-\mu_p)^{-1}$ and expand $Q^{-1}(\lambda)$ into a Taylor series in powers of $(\lambda-\mu_p)$, and next integrate over the circle $G(\mu_p)=\{\lambda:|\lambda-\mu_p|<r\}$, $0<r\ll1$, inside of which lies the only point $\mu_p$, then we obtain the following $N$ equations:
$$({\rm iii})\,\frac12R_q(\zeta_0,x)\ddot{Q}^{-1}(\mu_p)=R_p(\zeta_0,x)\frac{\zeta_1e^{i(\zeta_0-\zeta_2)n_+\mu_px}}{(\mu_p+\mu_p\zeta_2)}+\sum\limits_{n\not=p}R_n(\zeta_0,x)\frac{\zeta_1e^{i(\zeta_0-\zeta_2)n_+\mu_nx}}{(\mu_p+\mu_n\zeta_2)^2}$$
\begin{equation}
+
\sum\widehat{R}_m(\zeta_0,x)\left[\frac1{(\mu_p-\nu_n)^2}+\frac{\zeta_2e^{i(\zeta_0-\zeta_1)n_+\nu_mx}}{(\mu_p+\nu_m\zeta_1)^2}\right]\label{eq2.52}
\end{equation}
$$+\frac1{2\pi i}\int\limits_0^\infty\frac{p_2(i\tau\zeta_1,x)v_1^+(\tau,x)}{Q(i\tau\zeta_1)}\frac{d\tau}{\tau+i\zeta_2\mu_p}
+\frac1{2\pi i}\int\limits_0^\infty\frac{p_1(i\tau\zeta_2,x)v_2^+(i\tau,x)}{Q(i\tau\zeta_2)}\frac{d\tau}{\tau+i\zeta_1\mu_p}$$
$$+\frac1{2\pi i}\int\limits_0^\infty\frac{p_3(i\tau,x)v_1^+(i\tau,x)+p_4(i\tau,x)v_4(i\tau,x)}{Q(i\tau)}\frac{d\tau}{\tau+i\mu_p}$$
($1\leq p\leq N$). Using analogous considerations for the point $\nu_p$, we obtain another $M$ equations,
$$({\rm iv})\quad\frac12\widehat{R}_q(\zeta_0,x)\ddot{Q}^{-1}(\nu_q)=\sum\limits_nR_n(\zeta_0,x)\left[\frac1{(\nu_q-\mu_n)^2}+\frac{\zeta_1e^{i(\zeta_0-\zeta_2)n_+\mu_nx}}{(\nu_q+\mu_n\zeta_2)^2}\right]$$
\begin{equation}
\begin{array}{ccc}
{+\displaystyle R_q(\zeta_0,x)\frac{\zeta_2e^{i(\zeta_0-\zeta_1)n_+\nu_qx}}{(\nu_q+\nu_q\zeta_1)^2}+\sum\limits_{m\not=q}R_m(\zeta_0,x)\left[\frac1{(\nu_q-\nu_m)^2}+\frac{\zeta_2e^{i(\zeta_0-\zeta_1)n_+\nu_mx}}{(\nu_q+\nu_m\zeta_1)^2}\right]}\\
{\displaystyle+\frac1{2\pi i}\int\limits_0^\infty\frac{p_2(i\tau\zeta_1,x)v_1^+(i\tau,x)}{Q(i\tau\zeta_1)}\frac{d\tau}{\tau+i\zeta_2\nu_q}+\frac1{2\pi i}\int\limits_0^\infty\frac{p_1(i\tau\zeta_2,x)v_2^+(i\tau,x)}{Q(i\tau\zeta_2)}\frac{d\tau}{\tau+i\zeta_1\nu_q}}
\end{array}\label{eq2.53}
\end{equation}
$$+\frac1{2\pi i}\int\limits_0^\infty\frac{p_3(i\tau,x)v_1^+(i\tau,x)+p_4(i\tau)v_2(i\tau,x)}{Q(i\tau)}\frac{d\tau}{\tau+i\nu_q}\quad(1\leq q\leq M).$$

\begin{conclusion}
We obtained the {\bf closed linear system of singular integral equations} \eqref{eq2.51} -- \eqref{eq2.53} relative to $v_1^+(it,x)$, $v_2^+(it,x)$ and $\{R_n(\zeta_0,x)\}_1^N$, $\{R_m(\zeta_0,x)\}_1^M$, independent parameters of which are $\{p_k(it,x)\}_1^4$ and number sets $\{\mu_n\}_1^N\in l_{\zeta_0}$. $\{\nu_m\}_1^m\in\widehat{l}_{\zeta_0}$. Note that \eqref{eq2.44}, \eqref{eq2.46} imply that functions $\{p_k(\lambda,x)\}_1^4$ are expressed via scattering coefficients $s_1(\lambda)$ and $s_2(\lambda)$ \eqref{eq2.32}.
\end{conclusion}

\begin{conclusion}
Knowing the solution to this ststem, using formula \eqref{eq2.49}, we find $F_Q(\lambda,x)$ \eqref{eq2.48}, and thus function $F(\lambda,x)$ \eqref{eq2.47} also. This function coincides with $\psi_0^+(\lambda,x)$ when $\lambda\in\Omega_0^-$. Using now Theorem \ref{t1.4} and \eqref{eq1.65}, we find function $m(x)$ on the right half-axis $x\in\mathbb{R}_+$.
\end{conclusion}

\section{Dual scattering problem for the incident waves coming from ``$-\infty$''}\label{s3}

{\bf 3.1} Analogously to \eqref{eq2.1}, we expand every Jost solution $\{v_k(\lambda,x)\}_0^2$ in other Jost solutions $\{u_k(\lambda,x)\}_0^2$,
\begin{equation}
v_k(\lambda,x)=\sum\limits_{l=0}^2\widetilde{t}_{k,l}(\lambda)u_k(\lambda,x)\quad(0\leq k\leq2),\label{eq3.1}
\end{equation}
or
\begin{equation}
v(\lambda,x)=\widetilde{T}(\lambda,x)u(\lambda,x)\label{eq3.2}
\end{equation}
where, as usual, $v(\lambda,x)\stackrel{\rm def}{=}\col[v_0(\lambda,x),v_1(\lambda,x),v_2(\lambda,x)]$, $u(\lambda,x)\stackrel{\rm def}{=}\col[u_0(\lambda,x),$ $u_1(\lambda,x),u_2(\lambda,)]$, and $\widetilde{T}(\lambda)$ is the dual transition matrix,
\begin{equation}
\widetilde{T}(\lambda)\stackrel{\rm def}{=}\left[
\begin{array}{ccc}
\widetilde{t}_{0,0}(\lambda)&\widetilde{t}_{0,1}(\lambda)&\widetilde{t}_{0,2}(\lambda)\\
\widetilde{t}_{1,0}(\lambda)&\widetilde{t}_{1,1}(\lambda)&\widetilde{t}_{1,2}(\lambda)\\
\widetilde{t}_{2,0}(\lambda)&\widetilde{t}_{2,1}(\lambda)&\widetilde{t}_{2,2}(\lambda)
\end{array}\right]\label{eq3.3}
\end{equation}
Equations \eqref{eq2.2} and \eqref{eq3.2} imply that matrix $\widetilde{T}(\lambda)$ \eqref{eq3.3} is the inverse of $T(\lambda)$ \eqref{eq2.3},
\begin{equation}
\widetilde{T}(\lambda)T(\lambda)=I.\label{eq3.4}
\end{equation}

\begin{remark}\label{r3.1}
Matrix $\widetilde{T}(\lambda)$ is calculated explicitly, since {\rm (ii)} \eqref{eq2.14} implies that
$$\widetilde{T}(\lambda)=k^{-2}JT^*(\overline{\lambda})J$$
where $J$ is given by \eqref{eq2.12} and $k=n_-/n_+$. Hence we find that
\begin{equation}
\widetilde{T}(\lambda)=k^{-2}\left[
\begin{array}{ccc}
t_{0,0}^+(\lambda)&\zeta_1t_{2,0}^+(\lambda)&\zeta_2t_{1,0}^+(\lambda)\\
\zeta_2t_{0,2}^+(\lambda)&t_{2,2}^+(\lambda)&\zeta_1t_{1,2}^+(\lambda)\\
\zeta_1t_{0,1}^+(\lambda)&\zeta_2t_{2,1}^+(\lambda)&t_{1,1}^+(\lambda)
\end{array}\right].\label{eq3.5}
\end{equation}
\end{remark}

\begin{remark}\label{r3.2}
It is sufficient to consider one of the equalities \eqref{eq3.1} (e.g., for $k=0$),
\begin{equation}
v_0(\lambda,x)=\widetilde{t}_{0,0}(\lambda)u_0(\lambda,x)+\widetilde{t}_{0,1}(\lambda)u_1(\lambda,x)+\widetilde{t}_{0,2}(\lambda)u_2(\lambda,x),\label{eq3.6}
\end{equation}
and other equalities in \eqref{eq3.1} follow from this one upon substituting $\lambda\rightarrow\lambda\zeta_1$ and $\lambda\rightarrow\lambda\zeta_2$, besides, functions $\widetilde{t}_{k,s}(\lambda)$ satisfy equalities similar to \eqref{eq2.5}.
\end{remark}

Analogously to \eqref{eq2.8}, the following representations hold for Wronskians $W_{k,s}(u,\lambda,x)\stackrel{\rm def}{=}W\{u_k(\lambda,x),u_s(\lambda,x)\}$:
\begin{equation}
\begin{array}{ccc}
W_{1,2}(u,\lambda,x)=\sqrt3z_-\zeta_0u_0^+(\lambda,x);\quad W_{0,1}(u,\lambda,x)=\sqrt3z_-\zeta_2u_1^+(\lambda,x);\\
 W_{2,0}(u,\lambda,x)=\sqrt3z_-\zeta_1u_2^+(\lambda,x).
\end{array}\label{eq3.7}
\end{equation}

\begin{picture}(200,200)
\put(0,100){\vector(1,0){200}}
\put(100,0){\vector(0,1){200}}
\put(100,100){\vector(1,2){50}}
\put(100,100){\vector(-1,2){50}}
\put(100,100){\vector(-3,1){100}}
\put(100,100){\vector(3,1){100}}
\put(71,145){$\times$}
\put(75,156){$\mu_n\zeta_1$}
\put(105,170){$\Omega_0$}
\put(122,150){$\circ$}
\put(126,145){$\nu_m\zeta_2$}
\qbezier(153,120)(100,140)(49,120)
\end{picture}

\hspace{20mm} Fig. 8

Since $\widetilde{t}_{0,0}(\lambda)=t_{00}^+(\lambda)/k^{-2}$, then, in view of Lemma \ref{l2.6}, function $\widetilde{t}_{00}(\lambda)$ is holomorphic in the sector $\Omega_0$ and has in it a finite number of zeros of multiplicity two of the form $\{\mu_n\zeta_1\}\in l_{\zeta_1}$ ($\mu_n>0$), $\{\nu_m\zeta_2\}\in\widehat{l}_{\zeta_2}$ ($\nu_m<0$).

Rewrite equality \eqref{eq3.6} as
\begin{equation}
\widetilde{r}_0(\lambda)u_0(\lambda,x)=u_0(\lambda,x)+\widetilde{s}_1(\lambda)u_1(\lambda,x)+\widetilde{s}_2(\lambda)u_2(\lambda,x)\label{eq3.8}
\end{equation}
where
\begin{equation}
\widetilde{r}_0(\lambda)=\frac1{\widetilde{t}_{0,0}(\lambda)};\quad\widetilde{s}_1(\lambda)=\frac{\widetilde{t}_{0,1}(\lambda)}{\widetilde{t}_{0,0}(\lambda)};\quad\widetilde{s}_2(\lambda)=\frac{\widetilde{t}_{0,2}(\lambda)}{\widetilde{t}_{0,0}(\lambda)}.
\label{eq3.9}
\end{equation}
And since, due to \eqref{eq3.5},
\begin{equation}
\widetilde{t}_{0,0}=\frac{t_{0,0}^+(\lambda)}{k^2};\quad\widetilde{t}_{0,1}(\lambda)=\zeta_1\frac{t_{2,0}^+(\lambda)}{k^2};\quad\widetilde{t}_{0,2}(\lambda)=\zeta_2\frac{t_{1,0}^+(\lambda)}{k^2},\label{eq3.10}
\end{equation}
then
\begin{equation}
\widetilde{r}_0(\lambda)=\frac{k^2}{t_{0,0}^+(\lambda)};\quad\widetilde{s}_1(\lambda)=\zeta_1\frac{t_{2,0}^+(\lambda)}{t_{0,0}^+(\lambda)};\quad\widetilde{s}_2(\lambda)=\zeta_2\frac{t_{1,0}^+(\lambda)}{t_{0,0}^+(\lambda)}.\label{eq3.11}
\end{equation}

\begin{remark}\label{r3.3}
For $x\rightarrow-\infty$, function $\widetilde{r}_0(\lambda)v_0(\lambda,x)$ \eqref{eq3.8} has an asymptotic
$$\widetilde{r}_0(\lambda)v_0(\lambda,x)\rightarrow e^{iz_-\zeta_0x}+\widetilde{s}_1(\lambda)e^{iz_-\zeta_1x}+\widetilde{s}_2e^{iz_-\zeta_2x}\quad(z_-=\lambda n_-),$$
and thus incident (from $-\infty$) wave $e^{iz_-\zeta_0}$ has the {\bf reflected scattered wave} $s_1(\lambda)e^{iz_-\zeta_1x}+s_2(\lambda)e^{iz_-\zeta_2x}$ where $\widetilde{s}_1(\lambda)$, $\widetilde{s}_2(\lambda)$ are the {\bf scattering coefficients} of wave $e^{iz_-\zeta_0x}$. Asymptotic behavior of $\widetilde{r}_0(\lambda)v_0(\lambda,x)$ at $+\infty$ is
$$\widetilde{r}_0(\lambda)v_0(\lambda,x)\rightarrow\widetilde{r}_0(\lambda)e^{iz_+\zeta_0x}\quad(z_+=\lambda n_+),$$
therefore, $\widetilde{r}_0(\lambda)$ is the {\bf coefficient of transformation} of incident wave $e^{iz_-\zeta_0x}$ into $e^{iz_+\zeta_0x}$.

Equations \eqref{eq3.4} and {\rm(ii)} \eqref{eq2.14} imply that for the matrix $\widetilde{T}(\lambda)$ the following equality holds:
\begin{equation}
\widetilde{T}(\lambda)J\widetilde{T}^*(\overline{\lambda})=k^{-2}J.\label{eq3.12}
\end{equation}
Upon equating the $(1,1)$-elements in this equality, we have
$$\widetilde{t}_{0,0}(\lambda)\widetilde{t}_{0,0}^+(\lambda)+\zeta_1\widetilde{t}_{1,0}(\lambda)\widetilde{t}_{2,0}^++\zeta_2\widetilde{t}_{2,0}(\lambda)\widetilde{t}_{1,0}^+(\lambda)=k^{-2},$$
or, in terms of $\widetilde{r}_0(\lambda)$, $\widetilde{s}_1(\lambda)$, $\widetilde{s}_2(\lambda)$ \eqref{eq3.11},
\begin{equation}
1+\zeta_1\widetilde{s}_1(\lambda)\widetilde{s}_2^+(\lambda)+\zeta_2\widetilde{s}_2(\lambda)s_1^+(\lambda)=k^{-2}\widetilde{r}_0(\lambda)r_0^+(\lambda)\label{eq3.13}
\end{equation}
Equality \eqref{eq3.13} is the {\bf unitarity property of the dual scattering problem}.
\end{remark}

Conservation law of the direct scattering problem (Remark \ref{r2.3}) is given by
\begin{equation}
k^2r_0(\lambda)r_0^+(\lambda)=1+\zeta_1s_1(\lambda)s_2^+(\lambda)+\zeta_2s_2(\lambda)s_1^+(\lambda)\label{eq3.14}
\end{equation}
Since $\widetilde{r}_0(\lambda)=k^2r_0^+(\lambda)$, then $\widetilde{r}_0\widetilde{r}_0^+(\lambda)=k^4r_0(\lambda)r_0^+(\lambda)$, therefore, equality \eqref{eq3.13} becomes
\begin{equation}
1+\zeta_1\widetilde{s}_1(\lambda)\widetilde{s}_2^+(\lambda)+\zeta_2\widetilde{s}_2(\lambda)s_1^+(\lambda)=k^2r_0(\lambda)r_0^+(\lambda).\label{eq3.15}
\end{equation}
Equations \eqref{eq3.14}, \eqref{eq3.15} imply
\begin{equation}
\zeta_1s_1(\lambda)s_2^+(\lambda)+\zeta_2s_2(\lambda)s_1^+(\lambda)=\zeta_1\widetilde{s}_1(\lambda)s_2^+(\lambda)+\zeta_2\widetilde{s}_2(\lambda)s_1^+(\lambda).\label{eq3.16}
\end{equation}
Equality \eqref{eq3.16} one should consider as a {\bf property of reciprocity of direct and dual scattering problems}.
\vspace{5mm}

{\bf 3.2} Upon substituting $\lambda\rightarrow\lambda\zeta_1$ and $\lambda\rightarrow\lambda\zeta_2$ in \eqref{eq3.8}, we obtain
\begin{equation}
\begin{array}{ccc}
({\rm i})\quad\widetilde{r}_1(\lambda)v_1(\lambda,x)=\widetilde{s}_2(\lambda\zeta_1)u_0(\lambda,x)+u_1(\lambda,x)+s_1(\lambda\zeta_1)u_2(\lambda,x);\\
({\rm ii})\quad\widetilde{r}_2(\lambda)v_2(\lambda,x)=s_1(\lambda\zeta_2)u_0(\lambda,x)+s_2(\lambda\zeta_2)u_1(\lambda,x)+u_2(\lambda,x).
\end{array}\label{eq3.17}
\end{equation}
Relation (i) \eqref{eq3.17} corresponds to the incident (from $-\infty$) wave $e^{iz_-\zeta_1x}$, and (ii) \eqref{eq3.17}, correspondingly, to the wave $e^{iz_-\zeta_2x}$.

Using the form \eqref{eq3.8} of the function $v_0(\lambda,x)$, we calculate the Wronskian $W\{v_0(\lambda,x),u_k(\lambda,x)\}$ ($k=1$, 2), then, taking into account \eqref{eq3.7}, we have
\begin{equation}
\begin{array}{ccc}
({\rm i})\quad\widetilde{r}_0(\lambda)g_{0,1}(\lambda,x)=\sqrt3z_-\zeta_0u_1^+(\lambda,x)-\widetilde{s}_2(\lambda)\sqrt3z_-\zeta_0u_0^+(\lambda,x);\\
({\rm ii})\quad\widetilde{r}_0(\lambda)g_{0,2}(\lambda,x)=-\sqrt3z_-\zeta_1u_2^+(\lambda,x)+\widetilde{s}_1(\lambda)\sqrt3z_0\zeta_0u_0^+(\lambda,x),
\end{array}\label{eq3.18}
\end{equation}
where
\begin{equation}
g_{0,k}(\lambda,x)\stackrel{\rm def}{=}W\{v_0(\lambda,x),u_k(\lambda,x)\}\quad(k=1,2).\label{eq3.19}
\end{equation}

\begin{picture}(200,200)
\put(0,100){\vector(1,0){200}}
\put(100,0){\vector(0,1){200}}
\put(200,67){\line(-3,1){200}}
\put(0,67){\line(3,1){200}}
\put(120,150){$g_{0,2}(\lambda,x)$}
\put(45,150){$g_{0,1}(\lambda,x)$}
\qbezier(153,120)(100,140)(49,120)
\qbezier(100,23)(155,35)(168,125)
\qbezier(45,120)(30,70)(100,25)
\put(0,55){$u_1^+(\lambda,x)$}
\put(180,80){$u_2^+(\lambda,x)$}
\end{picture}

\hspace{20mm} Fig. 9

Functions $u_1^+(\lambda,x)$ and $u_2^+(\lambda,x)$ are holomorphic at the sectors $\Omega_2$ and $\Omega_1$, and functions $g_{0,k}(\lambda,x)$, correspondingly, at the sectors $\Omega_0\cap\Omega_k^-$ ($k=1$, 2). Using (i) \eqref{eq3.17}, calculate the Wronskian $W\{v_1(\lambda,x),u_k(\lambda,x)\}$ ($k=0$, 2), then taking into account \eqref{eq3.7}, we obtain
\begin{equation}
\begin{array}{ccc}
({\rm i})\quad\widetilde{r}_1(\lambda)g_{1,0}(\lambda,x)=-\sqrt3z_-\zeta_2u_1^+(\lambda,x)+\widetilde{s}_1(\lambda\zeta_1)\sqrt3z_-\zeta_2u_2^+(\lambda,x);\\
({\rm ii})\quad\widetilde{r}_1(\lambda)g_{1,2}(\lambda,x)=-s_2(\lambda\zeta_1)\sqrt3z_-\zeta_1u_2^+(\lambda,x)+\sqrt3z_-\zeta_0u_0^+(\lambda,x),
\end{array}\label{eq3.20}
\end{equation}
where
\begin{equation}
g_{1,k}(\lambda,x)=W\{v_1(\lambda,x),u_k(\lambda,x)\}\quad(k=0,2).\label{eq3.21}
\end{equation}

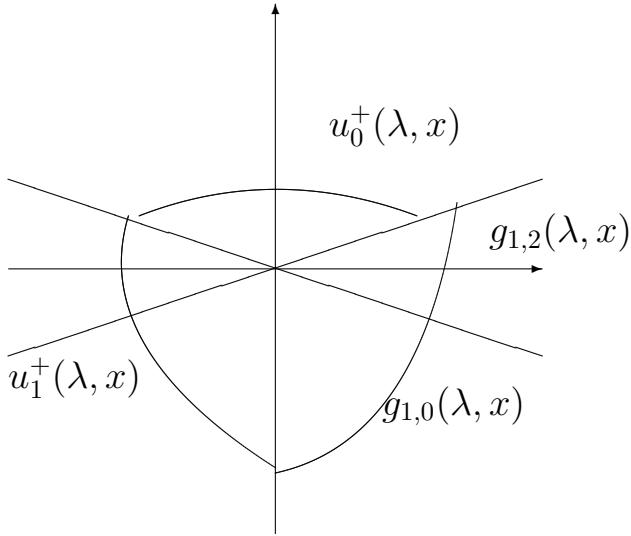
\begin{picture}(200,200)
\put(0,100){\vector(1,0){200}}
\put(100,0){\vector(0,1){200}}
\put(200,67){\line(-3,1){200}}
\put(0,67){\line(3,1){200}}
\put(120,150){$u_0^+(\lambda,x)$}
\qbezier(153,120)(100,140)(49,120)
\qbezier(100,23)(155,35)(168,125)
\qbezier(45,120)(30,70)(100,25)
\put(0,55){$u_1^+(\lambda,x)$}
\put(140,45){$g_{1,0}(\lambda,x)$}
\put(180,110){$g_{1,2}(\lambda,x)$}
\end{picture}

\hspace{20mm} Fig. 10

Functions $u_0^+(\lambda,x)$, $u_1^+(\lambda,x)$ are holomorphic at the sectors $\Omega_0$, $\Omega_2$, and functions $g_{1,k}(\lambda,x)$ are analytic, correspondingly, at the sectors $\Omega_1\cap\Omega_k^-$ ($k=0$, 2). Finally, using (ii) \eqref{eq3.17} and calculating $W\{v_2(\lambda,x),u_k(\lambda,x)\}$ ($k=0$, 1), in view of \eqref{eq3.7}, we have
\begin{equation}
\begin{array}{lll}
({\rm i})\quad\widetilde{r}_2(\lambda)g_{2,0}(\lambda,x)=-\widetilde{s}_2(\lambda\zeta_3)\sqrt3z_-\zeta_2u_1^+(\lambda,x)+\sqrt3z_-\zeta_1u_2^+(\lambda,x);\\
({\rm ii)}\quad\widetilde{r}_2(\lambda)g_{2,1}(\lambda,x)=\widetilde{s}_1(\lambda\zeta_2)\sqrt3z_-\zeta_2u_1^+(\lambda,x)-\sqrt3z_-\zeta_0u_0^+(\lambda,x),
\end{array}\label{eq3.22}
\end{equation}
where
\begin{equation}
g_{2,k}(\lambda,x)=W\{v_2(\lambda,x),u_k(\lambda,x)\}\quad(k=0,1).\label{eq3.23}
\end{equation}

\begin{picture}(200,200)
\put(0,100){\vector(1,0){200}}
\put(100,0){\vector(0,1){200}}
\put(200,67){\line(-3,1){200}}
\put(0,67){\line(3,1){200}}
\put(120,150){$u_0^+(\lambda,x)$}
\qbezier(153,120)(100,140)(49,120)
\qbezier(100,23)(155,35)(168,125)
\qbezier(60,115)(40,90)(60,85)
\qbezier(44,80)(40,60)(100,25)
\put(0,55){$g_{2,0}(\lambda,x)$}
\put(0,105){$g_{2,1}(\lambda,x)$}
\put(140,45){$u_2^+(\lambda,x)$}
\end{picture}

\hspace{20mm} Fig. 11

Functions $u_0^+(\lambda,x)$, $u_2^+(\lambda,x)$ are holomorphic at the sectors $\Omega_0$, $\Omega_1$, and $g_{2,k}(\lambda,x)$, correspondingly, at the sectors $\Omega_2\cap\Omega_k^-$ ($k=0$, 1). Using \eqref{eq3.19}, \eqref{eq3.21}, \eqref{eq3.23}, define holomorphic in corresponding sectors functions
$$\varphi_0^+(\lambda,x)\stackrel{\rm def}{=}u_0^+(\lambda,x)e^{iz_-\zeta_0x}\quad(\lambda\in\Omega_0);$$
$$\varphi_{1,2}(\lambda,x)\stackrel{\rm def}{=}\widetilde{r}_1(\lambda)g_{1,2}(\lambda,x)e^{iz_-\zeta_0x}\quad(\lambda\in\Omega_1\cap\Omega_2^-);$$
\begin{equation}
\varphi_{2,1}(\lambda,x)\stackrel{\rm def}{=}r_2(\lambda)g_{1,2}(\lambda,x)e^{iz_-\zeta_0x}\quad(\lambda\in\Omega_2\cap\Omega_1^-);\label{eq3.24}
\end{equation}
$$\varphi_{1,0}(\lambda,x)\stackrel{\rm def}{=}\widetilde{r}_1(\lambda)g_{1,0}(\lambda,x)e^{iz_-\zeta_2x}\quad(\lambda\in\Omega_1\cap\Omega_0^-);$$
$$\varphi_{2,0}(\lambda,x)\stackrel{\rm def}{=}r_2(\lambda)g_{2,0}(\lambda,x)e^{iz_-\zeta_1x}\quad(\lambda\in\Omega_2\cap\Omega_0^-).$$
Note that $\varphi_0^+(\lambda,x)$ is analytic at the sector $\Omega_0$, and functions $\varphi_{1,2}(\lambda,x)$, $\varphi_{2,1}(\lambda,x)$, $\varphi_{1,0}(\lambda,x)$, $\varphi_2(\lambda,0)$ are meromorphic at the corresponding sectors and have poles at the points $\{\mu_n\}$, $\{\nu_m\}$, $\{\nu_m\zeta_1\}$, $\{\mu_n\zeta_2\}$ correspondingly.

\begin{picture}(200,200)
\put(0,100){\vector(1,0){200}}
\put(100,0){\vector(0,1){200}}
\put(200,67){\line(-3,1){200}}
\put(0,67){\line(3,1){200}}
\put(120,150){$\varphi_0^+(\lambda,x)$}
\put(170,110){$\varphi_{1,2}(\lambda,x)$}
\qbezier(153,120)(100,140)(49,120)
\qbezier(100,23)(155,35)(168,80)
\qbezier(60,115)(40,90)(60,85)
\qbezier(150,120)(168,110)(168,80)
\qbezier(44,80)(40,60)(100,25)
\put(0,52){$\varphi_{2,0}(\lambda,x)$}
\put(0,105){$\varphi_{2,1}(\lambda,x)$}
\put(140,30){$\varphi_{1,0}(\lambda,x)$}
\end{picture}

\hspace{20mm} Fig. 12

Equations (ii) \eqref{eq3.20} and (ii) \eqref{eq3.22} imply two jump problems on the rays $i\widehat{l}_{\zeta_1}$ and $i\widehat{l}_{\zeta_2}$:
\begin{equation}
\begin{array}{ccc}
\sqrt3z_-\zeta_0\varphi_0^+(\lambda,x)-\varphi_{1,2}(\lambda,x)=\widetilde{p}_2(\lambda,x)u_2^+(\lambda,x)\quad(i\widehat{l}_{\zeta_1});\\
\sqrt3z_+\zeta_0\varphi_0^+(\lambda,x)+\varphi_{2,1}(\lambda,x)=p_1(\lambda,x)u_1^+(\lambda,x)\quad(il_{\zeta_2}),
\end{array}\label{eq3.25}
\end{equation}
where
\begin{equation}
\widetilde{p}_2(\lambda,x)\stackrel{\rm def}{=}\widetilde{s}_2(\lambda\zeta_1)\sqrt3z_-\zeta_1e^{iz_-\zeta_0x};\quad\widetilde{p}_1(\lambda,x)\stackrel{\rm def}{=}\widetilde{s}_1(\lambda\zeta_2)\sqrt3z_-\zeta_2e^{iz_-\zeta_0x}.\label{eq3.26}
\end{equation}
Boundary values of the function $\varphi_{1,2}(\lambda,x)$ on the ray $il_{\zeta_2}$ coincide with boundary values of $\varphi_{2,0}(\lambda,x)$ on the ray $il_{\zeta_1}$ since
$$\left.\varphi_{1,2}(\lambda,x)\right|_{\lambda=i\tau\zeta_2\in il_{\zeta_2}}=\frac1{r_0(i\tau)}g_{0,1}(\lambda,x)e^{-n_-\tau\zeta_2x};$$
$$\left.\varphi_{2,0}(\lambda,x)\right|_{\lambda=i\tau\zeta_1\in il_{\zeta_1}}=\frac1{r_0(i\tau)}g_{0,1}(\lambda,x)e^{-n_-\tau\zeta_2x}.$$
Similarly, it is determined that $\left.\varphi_{2,1}(\lambda,x)\right|_{\lambda\in il_{\zeta_1}}=\left.\varphi_{1,0}(\lambda,x)\right|_{\lambda\in il_{\zeta_2}}$. Therefore, defining $\widehat{\varphi}_{2,0}(\lambda,x)$ at the sector $\Omega_1\cap\Omega_0^-$ by the formula $\widehat{\varphi}_{2,0}(\lambda,x)\stackrel{\rm def}{=}\varphi_{2,0}(\widehat{\lambda},x)$, where $\widehat{\lambda}$ is the point symmetric relative to the axis $i\zeta_0$, we obtain an analytic continuation of the function $\varphi_{1,2}(\lambda,x)$ into the sector $\Omega_1$. Similarly defining $\widehat{\varphi}_{0,1}(\lambda,x)$ using symmetry at the sector $\Omega_2\cap\Omega_0^-$, we obtain a holomorphic continuation of the function $\varphi_{2,1}(\lambda,x)$ into the sector $\Omega_2$.

\begin{picture}(200,200)
\put(0,100){\vector(1,0){200}}
\put(100,0){\vector(0,1){200}}
\put(200,67){\line(-3,1){200}}
\put(0,67){\line(3,1){200}}
\put(120,150){$\varphi_0^+(\lambda,x)$}
\put(175,140){$i\widehat{l}_{\zeta_1}$}
\put(175,110){$\varphi_{1,2}(\lambda,x)$}
\qbezier(153,120)(100,140)(49,120)
\qbezier(100,23)(155,35)(168,80)
\qbezier(60,115)(40,90)(60,85)
\qbezier(150,120)(168,110)(168,80)
\qbezier(44,80)(40,60)(100,25)
\put(0,52){$\widehat{\varphi}_{0,1}(\lambda,x)$}
\put(0,105){$\widehat{\varphi}_{2,1}(\lambda,x)$}
\put(0,135){$i\widehat{l}_{\zeta_2}$}
\put(140,30){$\widehat{\varphi}_{0,2}(\lambda,x)$}
\end{picture}

\hspace{20mm} Fig. 13

As a result, we obtain holomorphic at the sectors $\Omega_0$, $\Omega_1$, $\Omega_2$ functions which at the sector $\Omega_1$ have poles of multiplicity 2 at the points $\{\mu_n\}$, $\{-\mu_n\zeta_1\}$ and at the sector $\Omega_2$, correspondingly, at $\{\nu_m\}$, $\{-\nu_m\zeta_2\}$. Multiply equality (i) \eqref{eq3.22} by $e^{iz_-\zeta_1x}$ and equality (i) \eqref{eq3.20} by $e^{iz_-\zeta_2x}$, and subtract the obtained, then we have
\begin{equation}
\varphi_{2,0}(\lambda,x)-\varphi_{1,0}(\lambda,x)=\widetilde{p}_3(\lambda,x)u_1^+(\lambda,x)+\widetilde{p}_4(\lambda,x)u_2^+(\lambda,x)\label{eq3.26}
\end{equation}
where
\begin{equation}
\begin{array}{ccc}
\widetilde{p}_3(\lambda,x)=\sqrt3z_-\zeta_2\left[\widetilde{s}_2(\lambda\zeta_2)e^{iz_-\zeta_1x}-e^{iz_-\zeta_2x}\right];\\
\widetilde{p}_4(\lambda,x)=\sqrt3z_-\zeta_1\left[e^{iz_-\zeta_1x}-\widetilde{s}_1(\lambda\zeta_1)e^{iz_-\zeta_2x}\right].
\end{array}\label{eq3.27}
\end{equation}
Define the function
\begin{equation}
\Phi(\lambda,x)=\left\{
\begin{array}{lllll}
\varphi_0^+(\lambda,x)\quad(\lambda\in\Omega_0);\\
\varphi_{2,1}(\lambda,x)\quad(\lambda\in\Omega_2\cap\Omega_1^-);\\
\widehat{\varphi}_{0,1}(\lambda,x)\quad(\lambda\in\Omega_2\cap\Omega_0^-);\\
\widehat{\varphi}_{0,2}(\lambda,x)\quad(\lambda\in\Omega_1\cap\Omega_0^-);\\
\varphi_{1,2}(\lambda,x)\quad(\lambda\in\Omega_1\cap\Omega_2^-).
\end{array}\right.\label{eq3.28}
\end{equation}
and let
\begin{equation}
{\it\Phi}_Q(\lambda,x)=\Phi(\lambda,x)/Q(\lambda)\label{eq3.29}
\end{equation}
where $Q(\lambda)$ is an entire function not having zeros, which for $\lambda\rightarrow\infty$ grows faster than $\exp\{|\lambda|(\sigma_+(0)+\sigma_-(0)\}$, e.g., $Q(\lambda)=\exp\lambda^2(\sigma_+(0)+\sigma_-(0))$. Function ${\it\Phi}_Q(\lambda,x)$ is holomporphic at the sector $\Omega_0$ and meromorphic at the sectors $\Omega_1$, $\Omega_2$, and vanishes as $\lambda\rightarrow0$, therefore \cite{16,17}, ${\it\Phi}_Q(\lambda,x)$ \eqref{eq3.29} is determined by its jumps on the rays $i\widehat{l}_{\zeta_1}$, $i\widehat{l}_{\zeta_2}$, $i\widehat{l}_{\zeta_0}$ via coefficients before the second-order poles at the points $\{\mu_n\}$, $\{-\mu_n\zeta_1\}$, $\{\nu_m\}$, $\{-\nu_m\zeta_2\}$, and thus
$${\it\Phi}_Q(\lambda,x)=\sum\limits_n\frac{R'_n(\zeta_0,x)}{(\lambda-\mu_n)^2}+\sum\limits_n\frac{R'_n(\zeta_1,x)}{(\lambda+\mu_n\zeta_1)^2}+\sum\limits_m\frac{\widehat{R}'_m(\zeta_0,x)}{(\lambda-\nu_m)^2}+\sum\limits_m
\frac{\widehat{R}'_m(\zeta_2,x)}{(\lambda+\nu_m\zeta_2)^2}$$
\begin{equation}
+\frac1{2\pi i}\int\limits_0^\infty\frac{\widetilde{p}_1(-i\tau\zeta_2,x)u_2^+(-i\tau,x)}{Q(-i\tau\zeta_2)}\frac{d\tau}{\tau-i\zeta_1\lambda}+\frac1{2\pi i}\int\limits_0^\infty\frac{p_2(-i\tau\zeta_1,x)u_1^+(-i\tau,x)}{Q(-i\tau\zeta_1)}\frac{d\tau}{\tau-i\zeta_2\lambda}\label{eq3.30}
\end{equation}
$$+\frac1{2\pi i}\int\limits_0^\infty\frac{p_3(-i\tau,x)u_1^+(-i\tau,x)+p_4(-i\tau,x)u_2^+(-i\tau,x)}{Q(-i\tau)}\frac{d\tau}{\tau-i\lambda}.$$
Coefficients $R_n(\zeta_0,x)$ and $R'_n(\zeta_1,x)$ are dependent. Really, taking into account that $\left.g_{1,2}(\lambda,x)\right|_{\lambda=\lambda\zeta_2}=g_{0,1}(\lambda,x)$ and using \eqref{eq3.24}, it is easy to show that
\begin{equation}
R'_n(\zeta_1,x)=\zeta_1e^{i(\zeta_0-\zeta_2)\mu_nnx}R'_n(\zeta_0,x).\label{eq3.31}
\end{equation}
Analogously, $\left.g_{2,1}(\lambda,x)\right|_{\lambda=\lambda\zeta_1}=g_{0,2}(\lambda,x)$ implies that
\begin{equation}
\widehat{R}'_m(\zeta_2,x)=\zeta_2e^{i(\zeta_0-\zeta_1)\nu_mn_-x}\widehat{R}'_m(\zeta_0,x).\label{3.32}
\end{equation}
Therefore, equality \eqref{eq3.30} becomes
$${\it\Phi}_Q(\lambda,x)=\sum\limits_nR'_n(\zeta_0,x)\left(\frac1{(\lambda-\mu_n)^2}+\frac{\zeta_1e^{i(\zeta_0-\zeta_2)\mu_nn_-x}}{(\lambda+\mu_n\zeta_1)^2}\right)$$
$$+\sum\limits_m\widehat{R}_m(\zeta_0,x)\left(\frac1{(\lambda-\nu_m)^2}+
\frac{\zeta_2e^{i(\zeta_0-\zeta_1)n_-\nu_mx}}{(\lambda-\nu_m\zeta_2)^2}\right)$$
\begin{equation}
+\frac1{2\pi i}\int\limits_0^\infty\frac{p_1(-i\tau\zeta_2,x)u_2^+(-i\tau,x)}{Q(-i\tau\zeta_2)}\frac{d\tau}{\tau+i\zeta_1\lambda}+\frac1{2\pi i}\int\limits_0^\infty\frac{p_2(-i\tau\zeta_1,x)u_1^+(-i\tau,x)}{Q(-i\tau\zeta_1)}\frac{d\tau}{\tau+i\zeta_2\lambda}\label{eq3.33}
\end{equation}
$$+\frac1{2\pi i}\int\limits_0^\infty\frac{p_3(-i\tau,x)u_1^+(-i\tau,x)+p_4(-i\tau,x)u_2^+(-i\tau,x)}{Q(-i\tau)}\frac{d\tau}{\tau+i\lambda}.$$
For $\lambda\in\Omega_0$, ${\it\Phi}_Q(\lambda,x)=\varphi_0^+(\lambda,x)/B(\lambda)=u_0^+(\lambda,x)e^{iz_-\zeta_0x}/B(\lambda)$. Therefore, upon calculating boundary values for $\lambda\rightarrow-it\zeta_2\in i\widehat{l}_{\zeta_2}$ and $\lambda\rightarrow-it\zeta_1\in\widehat{l}_{\zeta_1}$ in \eqref{eq3.33} and using the Sokohtski formula \cite{16,17} for Cauchy type integrals, we obtain
$$({\rm i})\,u_1^+(-it,x)e^{t\zeta_2n_-tx}Q^{-1}(-it\zeta_2)=\sum\limits_nR'_n(\zeta_0,x)\left[\frac1{(it\zeta_2+\mu)^2}+\frac{\zeta_1e^{i(\zeta_0-\zeta_1)\mu_n}}{(it\zeta_2-\mu_n\zeta_1)^2}\right]$$
\begin{equation}
+\sum\limits_m\widehat{R}_m(\zeta_0,x)\left[\frac1{(it\zeta_2+\nu_m)^2}+\frac{\zeta_2e^{i(\zeta_0-\zeta_2)\nu_mx}}{(it\zeta_2+\nu_m\zeta_2)^2}+\frac12\frac{p_1(it\zeta_2,x)u_2^+(-it,x)}{Q(-it\zeta_2)}\right.\label{eq3.35'}
\end{equation}
$$+\frac1{2\pi i}/\hspace{-4.4mm}\int\limits_0^\infty\frac{p_1(-i\tau\zeta_2,x)u_2^+(-it,x)}{Q(-i\tau\zeta_2)}\frac{d\tau}{\tau-t}+\frac1{2\pi i}\int\limits_0^\infty\frac{p_1(-i\tau\zeta_1,x)u_1^+(-i\tau,x)}{Q(-i\tau\zeta_1)}\frac{d\tau}{\tau+\zeta_1t}$$
\begin{equation}
\begin{array}{ccc}
{\displaystyle\left.+\frac1{2\pi i}\int\limits_0^\infty\frac{p_3(-i\tau,x)u_1^+(-i\tau,x)+p_4(-i\tau,x)u_2^+(-i\tau,x)}{Q(-i\tau)}\frac{d\tau}{\tau-\zeta_2t}\right];}\\
{\displaystyle({\rm ii})\quad u_2^+(-it,x)e^{t\zeta_1n_-tx}=\sum\limits_nR'_n(\zeta_0,x)\left[\frac1{(it\zeta_1+\mu_n)^2}+\frac{\zeta_1e^{i(\zeta_0-\zeta_1)\mu_nx}}{(it\zeta_1-\mu_n\zeta_1)^2}\right]}
\end{array}\label{eq3.34}
\end{equation}
$$+\sum\limits_m\widehat{R}'_m(\zeta_0,x)\left[\frac1{(it\zeta_1+\nu_m)^2}+\frac{\zeta_2e^{i(\zeta_0-\zeta_2)\nu_mx}}{(it\zeta_1-\nu_m\zeta_2)^2}\right.$$
$$+\frac1{2\pi i}\int\limits_0^\infty\frac{p_1(-i\tau\zeta_2,x)u_2^+(-i\tau,x)}{Q(-i\tau\zeta_2)}\frac{d\tau}{\tau+\zeta_2t}$$
$$+\frac12\frac{p_2(-it\zeta_1,x)u_1^+(-it,x)}{Q(-it\zeta_1)}+\frac1{2\pi i}/\hspace{-4.4mm}\int\limits_0^{2\pi}\frac{p_2(-i\tau\zeta_1,x)u_1(-i\tau,x)}{Q(-i\tau\zeta_1)}\frac{d\tau}{\tau-t}$$
$$\left.+\frac1{2\pi i}\int\limits_0^\infty\frac{p_3(-i\tau,x)u_1^+(-i\tau,x)+p_4(-i\tau,x)u_2^+(-i\tau,x)}{Q(-i\tau)}\frac{d\tau}{\tau-\zeta_1t}\right].$$
Multiply equality \eqref{eq3.33} by $(\lambda-\mu_p)^{-1}$ and expand $Q^{-1}(\lambda)$ into a Taylor series in powers of $(\lambda-\mu_p)$, then upon integrating over a circle $C_r(\mu_p)=\{\lambda:|\lambda-\mu_p|=r\}$, $0<r\ll1$, inside of which there is only one point $\mu_p$, we obtain $N$ equations:
$$({\rm iii})\quad\frac12R'_p(\zeta_0,x)\ddot{Q}^{-1}(\mu_p)=R'_p(\zeta_0,x)\cdot\frac{\zeta_1e^{i(\zeta_0-\zeta_1)n_-\mu_px}}{(\mu_p+\mu_p\zeta_1)^2}+\sum\limits_{n\not=p}R'_n(\zeta_0,x)\left[\frac1{(\mu_p-\mu_n)^2}\right.$$
\begin{equation}
\begin{array}{ccc}
{\displaystyle\left.+\frac{\zeta_1e^{i(\zeta_0-\zeta_2)n_-\mu_nx}}{(\mu_p-\mu_n\zeta_1)^2}\right]+\sum\limits_mR'_m(\zeta_0,x)\left[\frac1{(\mu_p-\nu_m)^2}+\frac{\zeta_2e^{i(\zeta_0-\zeta_1)n_-\nu_mx}}{(\mu_p-\nu_m\zeta_2)^2}\right]}\\
{\displaystyle+\frac1{2\pi i}\int\limits_0^\infty\frac{p_1(-i\tau\zeta_2,x)u_2^+(-i\tau,x)}{Q(-i\tau\zeta_2)}\frac{d\tau}{\tau+i\zeta_1\mu_p}}\\
{\displaystyle+\frac1{2\pi i}\int\limits_0^\infty\frac{p_1(-i\tau\zeta_1,x)u_1^+(-i\tau,x)}{Q(-i\tau\zeta_1)}\frac{d\tau}{\tau+i\zeta_2\mu_p}}
\end{array}\label{eq3.35}
\end{equation}
$$+\frac1{2\pi i}\int\limits_0^\infty\frac{p_3(-i\tau,x)u_1^+(-i\tau,x)+p_4(-i\tau,x)u_2^+(-i\tau,x)}{Q(-i\tau)}\frac{d\tau}{\tau+i\mu_p}$$
($1\leq p\leq N$).
Applying analogous considerations, we obtain another $M$ equations:
$$({\rm iv})\quad\frac12\widehat{R}_q(\zeta_0,x)\ddot{Q}^{-1}(\nu_q)=\sum\limits_nR'_n(\zeta_0,x)\left[\frac1{(\nu_q-\mu_n)^2}+\frac{\zeta_1e^{i(\zeta_0-\zeta_2)n_-\mu_n}}{(\nu_q+\mu_n\zeta_1)^2}\right]$$
\begin{equation}
\begin{array}{ccc}
{\displaystyle+\widehat{R}'_q(\zeta_0,x)\frac{\zeta_2e^{i(\zeta_0-\zeta_1)n_-\nu_mx}}{(\nu_q-\nu_q\zeta_2)^2}+\sum\limits_{m\not=q}\widehat{R}'_m(\zeta_0,x)\left[\frac1{(\nu_q-\nu_m)^2}+\frac{\zeta_2e^{i(\zeta_0-\zeta_1)n_-\nu_mx}}{(\nu_q-
\nu_m\zeta_2)^2}\right]}\\
{\displaystyle+\frac1{2\pi i}\int\limits_0^\infty\frac{p_1(-i\tau\zeta_2,x)u_2^+(-i\tau,x)}{Q(-i\tau\zeta_2)}\frac{d\tau}{\tau+i\zeta_1\nu_q}+\frac1{2\pi i}\int\limits_0^\infty\frac{p_2(-i\tau\zeta_1,x)u_1^+(-i\tau,x)}{Q(-i\tau\zeta_1)}\frac{d\tau}{\tau+i\zeta_2\nu_q}}
\end{array}\label{eq3.36}
\end{equation}
$$+\frac1{2\pi i}\int\limits_0^\infty\frac{p_3(-i\tau,x)u_1^+(-i\tau,x)+p_4(-i\tau,x)u_2^+(-i\tau,x)}{Q(-i\tau)}\frac{d\tau}{\tau+i\nu_q},\quad1\leq q\leq M.$$

\begin{conclusion}
We obtained a {\bf closed system of linear singular integral equations} {\rm (i)}, {\rm (ii)} \eqref{eq3.34}, {\rm (iii)} \eqref{eq3.35}, {\rm (iv)} \eqref{eq3.36} relative to $u_1^+(-it,x)$, $u_2^+(-it,x),$ $\{R'_n(\zeta_0,x)\}_1^N$, $\{R_m(\zeta_0,x)\}_1^M$ with free parameters $\{p_k(\lambda,x)\}_1^4$ and sets of points $\{\mu_n\}_1^N$, $\{\nu_m\}_1^M$. Functions $\{p_k(\lambda,x)\}_1^4$ are expressed explicitly via scattering coefficients $\widetilde{s}_1(\lambda)$, $\widetilde{s}_2(\lambda)$ \eqref{eq3.5}.
\end{conclusion}

\begin{conclusion}
Knowing the solution $u_1^+(-it,x)$, $u_2^+(-it,x)$, $\{R'_n(\lambda,x)\}_1^N$, $\{\widehat{R}'_m(\zeta_0,x)\}_1^M$ to system \eqref{eq3.34} -- \eqref{eq3.36}, using formula \eqref{eq3.30}, we define function ${\it\Phi}_Q(\lambda,x)$, whence we find function $u_0^+(\lambda,x)$ in sector $\Omega_0$. Finally, from Theorem \ref{t1.6} using \eqref{eq1.73} we find function $m(x)$ on the left half-axis, $x\in\mathbb{R}_-$.
\end{conclusion}

\begin{remark}\label{r3.4}
Study of the basic systems of equations for direct \eqref{eq2.51} -- \eqref{eq2.53} and dual \eqref{eq3.34} -- \eqref{eq3.35} scattering problems has to be continued.
\end{remark}

\begin{remark}\label{r3.5}
Sets
\begin{equation}
\{s_1(\lambda),s_2(\lambda),\{\mu_n\}_1^N,\{\nu_m\}_1^M;\{\widetilde{s}_1(\lambda),\widetilde{s}_2(\lambda),\{\mu_n\}_1^N,\{\nu_m\}_1^M\}\label{eq3.37}
\end{equation}
are data of the direct and dual scattering problem, besides, $\{s_1(\lambda)\}_1^2$ and $\{\widetilde{s}_1(\lambda)\}_1^2$ satisfy relations of unitarity \eqref{eq2.33} and \eqref{eq3.13}, and relation of duality \eqref{eq3.16}. Description of the scattering data \eqref{eq3.37} of direct and dual scattering problem for a cubic string that has a form of step one has to study separately.
\end{remark}

\renewcommand{\refname}{ \rm 
\end{Large}

\centerline{\bf References}}
\begin{thebibliography}{99}

\bibitem{1}
V. Marchenko,  {\it Sturm-Liouville Operators and Applications} Ser. ``Operator Theory: Advances and Applications'', Vol. 22, Birkh$\ddot{\rm a}$user Basel, 1986, xi + 367 pp.
\bibitem{2}
M. G. Krein, ``Solution of the inverse Sturm -- Liouville problem'', {\it Dokl. Acad. Nauk SSSR (N. S.)}, 76 (1951), 21 -- 24 (in Russian).
\bibitem{3}
B. M. Levitan, {\it Inverse Sturm -- Liouville Problems}, Utrecht, The Netherlands, VNU Science Press, 1987, 240 pp.
\bibitem{4}
P. Deift, E. Trubovitz, ``Inverse scattering on the line'', {\it Pure Appl. Math.}, 32 (1979), 121 -- 251.
\bibitem{5}
S. Novikov, S. V. Manakov, L. P. Pitaevskii, V. E. Zakharov, {\it Theory of Solitons. The Inverse Scattering Method} Contemp. Soviet Math. Consultance Bureau [Plenum], New York, 1984.
\bibitem{6}
A. Constantin, R. I.  Ivanov, J. Lenells, ``Inverse scattering transform for the
Degasperis-Procesi equation'', {\it Nonlinearity}, 23 (2010), 2559 -- 2575.
\bibitem{7}
A. Degasperis, D. D. Holm, A. Hone, ``A new integrable equation with peakon solutions'', {\it Theoret. Math. Phys}, 133:2 (2002), 1463 -- 1474.
\bibitem{8}
A. Degasperis, M. Processi, ``Asymptotic Integrability''. In A. Degasperis et al., Eds. Symmetry and Perturbation Theory, World Sientific, Singapore (1999), 23 -- 37.
\bibitem{9}
A. Boutet de Monvel, D. Shepelsky, ``A Riemann -- Hilbert approach for the Degasperis -- Procesi equations'', {\it Nonlinearity}, 26 (2013), 2081 -- 2107.
\bibitem{10}
J. Kohlenberg, H. Lundmark, J. Szmigielski, ``The inverse spectral problem for the discrete cubic string'', {\it Inverse Problems}, 23 (2007), 99 -- 121.
\bibitem{11}
M. G. Krein, ``Determination of the density of an inhomogenous symmetric string from its frequency spectrum'', {\it Dokl. Acad. Nauk SSSR}, 76 (1951), 345 -- 348 (in Russian).
\bibitem{12}
I. S. Kac, M. G. Krein, ``On the spectral function of the string'', {\it Amer. Math. Soc. Transl.}, Ser. 2, 103 (1974), 19 -- 102.
\bibitem{13}
H. Lundmark, J. Smiegielski, ``Multi-peakon solution of the Degasperis -- Processi equations'', {\it Inverse Problems}, 19 (2003), 1241 -- 1245.
\bibitem{14}
V. A. Zolotarev, ``Inverse scattering problem for a third-order operator with local potential'', {\it Journal of Differential Equations}, 379 (2024), 207 -- 257.
\bibitem{15}
V. A. Zolotarev, ``Inverse scattering problem for third order differential operators on the whole axis'', {\it Journal of Mathematical Analysis and Applications}, 546:2 (2025), 129331.
\bibitem{16}
F. D. Gahov, {\it Boundary Value Problems}, Elsevier, 2014, 584 pp.
\bibitem{17}
N.I. Muskhelishvili, {\it Singular Integral Equations}, Springer, Dordrecht, 1958, xiv + 441 pp.

\end{thebibliography}
\end{document}